\documentclass{amsart}
\usepackage[latin1]{inputenc}
\usepackage{amssymb}
\usepackage{amsmath}
\usepackage{enumerate}
\usepackage{epsfig,color,xcolor}

\usepackage[all]{xy}
\usepackage{pgf,tikz}
\usepackage{mathrsfs}
\usetikzlibrary{arrows}

\newtheorem{theorem}{Theorem}[section]

\newtheorem{prop}[theorem]{Proposition}
\newtheorem{proposition}[theorem]{Proposition}
\newtheorem{corollary}[theorem]{Corollary}

\newtheorem{definition}[theorem]{Definition}

\newtheorem{rem}[theorem]{Remark}
\newtheorem{remark}[theorem]{Remark}

\numberwithin{equation}{section}
\newcommand{\rk}{{\rm rank}}
\newcommand{\id}{\mbox{id}}

\newcommand{\MW}{\mbox{MW}}

\newcommand{\lra}{\leftrightarrow}

\newcommand{\ra}{\rightarrow}

\newcommand{\C }{ \mathbb{C}}

\newcommand{\Z}{\mathbb{Z}}
\newcommand{\Q}{\mathbb{Q}}

\newcommand{\Aut}{\mbox{Aut}}

\makeatletter
\def\blfootnote{\xdef\@thefnmark{}\@footnotetext}

\title[Generalized Shioda--Inose structures of order 3]{Generalized Shioda--Inose structures of order 3}
\begin{document}
	
	\author{Alice Garbagnati}
	\address{Alice Garbagnati, Dipartimento di Matematica, Universit\`a Statale di Milano,
		via Saldini 50, I-20133 Milano, Italy}
	\email{alice.garbagnati@unimi.it}
	
	\author{Yulieth Prieto Monta\~{n}ez}
	\address{Yulieth Prieto Monta\~{n}ez, The Abdus Salam International Centre for Theoretical Physics, Str. Costiera, 11, 34151 Trieste TS, Italy}
	\email{yprieto@ictp.it}


	\begin{abstract}
		A Shioda--Inose structure is a geometric construction which associates to an Abelian surface a projective K3 surface in such a way that their transcendental lattices are isometric. This geometric construction was described by Morrison by considering special symplectic involutions on the K3 surfaces. After Morrison several authors provided explicit examples. The aim of this paper is to generalize Morrison's results and some of the known examples to an analogous geometric construction involving not involutions, but order 3 automorphisms. Therefore we  define the generalized Shioda--Inose structures of order 3, we identify the K3 surfaces and the Abelian surfaces which appear in these structures and we provide explicit examples. \end{abstract}	
	\maketitle

\section*{Introduction}

The surfaces with trivial canonical bundle are the K3 surfaces and the 2-dimensional tori.
Their geometry presents several common aspects, and many of the important results on K3 surfaces have analogues for 2-dimensional tori.

The most famous relation between K3 surfaces and 2-dimensional tori is given by the Kummer construction: each 2-dimensional torus $A$ is endowed with an involution $\iota_A$, which fixes 16 points. The quotient $A/\iota_A$ is a singular surface, blowing up its singular points one finds a K3 surface, classically called the Kummer surface of $A$ and denoted in the following by $Km_2(A)$. The Kummer construction associates to each 2-dimension torus a K3 surface, but not to each K3 surface a 2-dimensional torus, indeed the family of the K3 surfaces is 20 dimensional, the one of the Kummer surfaces is 4 dimensional. 

The Kummer construction is interesting also because it creates a relation between the Hodge structures of the transcendental lattices of the torus and of its Kummer surface: the nowhere vanishing 2-holomorphic form on the torus $A$ induces the nowhere vanishing 2-holomorphic form on $Km_2(A)$ and hence the Hodge structure on the transcendental lattice $T_A$ induces the one on $T_{Km_2(A)}$. Moreover $T_{Km_2(A)}\simeq T_A(2)$, i.e. the bilinear form of $T_{Km_2(A)}$ is the one of $T_A$ multiplied by 2.

In \cite{SI}, Shioda and Inose relate a few specific 2-dimensional tori $A$ with K3 surfaces $X$ such that $T_A\simeq T_X$. The disadvantage with respect to the Kummer construction is the fact that Shioda and Inose consider only particular tori, but the advantage is the stronger relation between the transcendental lattices of the surfaces.

In \cite{Mo}, Morrison extends the results of Shioda and Inose to all the projective 2-dimensional tori. Morrison defines the ``Shioda--Inose structures", see Definition \ref{def: SI order 2}, which are the data of an Abelian surface $A$, a projective K3 surface $X$ and an involution $\iota_X$ on $X$ such that $T_A\simeq T_X$ and $X/\iota_X$ is birational to $Km_2(A)$; then he proves that for all the Abelian surfaces $A$ there exists a Shioda--Inose structure. Hence every projective Kummer surface is rationally $2:1$ covered both by an Abelian surface $A$ and by a K3 surface $X$ such that $T_A\simeq T_X$. The K3 surfaces $X$ that are members of a Shioda--Inose structure can be identified by properties of their transcendental lattices and N\'eron--Severi groups: a projective K3 surface is a part of a Shioda--Inose structure if and only if its N\'eron--Severi group primitively contains the lattice $E_8(-1)\oplus E_8(-1)$. Moreover, an involution $\iota_X$ is the involution of a Shioda--Inose structure if and only if the induced isometry $\iota_X^*$ acts on the N\'eron--Severi group switching these two orthogonal copies of $E_8(-1)$. By the classic construction of the moduli space of the lattice polarized K3 surfaces, the families of projective K3 surfaces admitting a Shioda--Inose structure have dimension 3 (as have the families of the Abelian surfaces and of the Kummer surfaces).

After Morrison, several authors describe explicitly K3 surfaces $X$ appearing in a Shioda--Inose structure, see e.g.  \cite{CD}, \cite{GL}, \cite{vGS}, \cite{K}, \cite{Sc} where 3-dimensional families of these type of K3 surfaces are described. In particular, in the latter three papers the involutions $\iota_X$ of the Shioda-Inose structures are constructed as a translation by a 2-torsion section on specific elliptic fibration.\\

As we saw, the first example of relations between 2-dimensional tori and K3 surfaces is provided by the Kummer construction, so it is obtained by considering an order 2 automorphism of the 2-dimensional tori, such that the minimal resolution of the quotient is a K3 surface. The Shioda--Inose structures introduced by Morrison have the Kummer construction as starting point. In order to generalize these structures one has first to consider generalizations of the Kummer construction. In \cite{Fu} Fujiki classifies the finite groups $G$ of automorphisms which act on a 2-dimensional torus in such a way that the minimal resolution of $A/G$ is a K3 surface. In particular, the cyclic ones have order 2,3,4,6 and we will denote by $Km_n(A)$ the K3 surface obtained as desingularization of a torus $A$ by the action of the cyclic group $\Z/n\Z$. In \cite{Be} a lattice theoretic description of the K3 surfaces $Km_n(A)$ is provided for each $n$ and in \cite{Ba} a classification of the families of the projective K3 surfaces of type $Km_3(A)$ is given. This essentially generalizes the known results for the Kummer construction of $Km_2(A)$. The main difference with respect to the case of order 2 is that not all the 2-dimensional tori admit automorphisms of order $n$ as required if $n>2$, so this construction allows to associate a K3 surface $Km_n(A)$ only to certain specific 2-dimensional tori.

To generalize the Shioda--Inose construction to higher order, one needs to find a K3 surface $X$ which is birational to a $\Z/n\Z$-Galois cover of $Km_n(A)$ and such that $T_X\simeq T_A$. In particular, this implies that both $A$ and $X$ admit an order $n$ automorphism of a certain prescribed type.
In \cite{CO} a generalization of the Shioda--Inose structures was proposed, but without the condition on the transcendental lattices. In this paper we discuss a stronger generalization, which includes the condition that the transcendental lattice of the K3 surface is isometric to the one of the Abelian surface; without this condition, even in the case of order 2, there are many different families of K3 surfaces which can be related with the same family of Abelian surfaces, which is not the case in the original Shioda--Inose construction proposed by Morrison and by Shioda and Inose. So we are interested in requiring a relation between the transcendental lattices of the Abelian surface and of the K3 surface in the Shioda--Inose structure.

We prove results analogous to the ones proved by Morrison for generalized Shioda--Inose structure of order 3 (see Theorem \ref{Main theroem intro}) and we also briefly discuss the difficulties in generalizing the notion of Shioda--Inose structure to the case of order 6 (see Section \ref{subsubsection: remarks n=6} and Remark \ref{rem: quotients different groups isomorphic}). In particular we prove that the generalized Shioda--Inose structures of order 3 have essentially all the good properties of the Shioda--Inose structures of order 2, but we also prove that it is not possible to generalize completely the Shioda--Inose construction to the order 6 case. For the order 6 case we prove that there exist an Abelian surface $A$ and a K3 surface $X$ with isometric transcendental lattices, such that $A$ admits an order 6 automorphism $\tau$ with $A/\tau$ birational to $Km_6(A)$, but $X$ does not admit an automorphism of order 6. Hence one can not construct a generalization of the Shioda--Inose structure for the case of order 6 requiring that the transcendental lattices of the Abelian surface and the of K3 surface are isometric. It is possible that there exist generalizations of the Shioda--Inose structures for the order 6 case which require a slightly different relation between the transcendental lattices of the involved surfaces. \\

In Section \ref{subsec: Kummer} we recall the main and well known results about the construction of Kummer surfaces $Km_2(A)$ and generalized Kummer surfaces of order 3 $Km_3(A)$. To define a (generalized) Shioda--Inose structure of order 3 on a K3 surface $X$ one has to consider certain order 3 automorphisms $\sigma_X$ on $X$, which have the property that the minimal resolution of $X/\sigma_X$ is a K3 surface, isomorphic to $Km_3(X)$. The condition that the minimal resolution of $X/\sigma_X$ is a K3 surface forces the automorphism $\sigma_X$ to be symplectic and in Section \ref{subsec: aymplectic autom} we recall the main results on symplectic automorphisms of order 2 and 3 on K3 surfaces. 
In Section \ref{subsec: Shioda Inose} we first recall the main results on the Shioda--Inose structures (of order 2) which were already proved in \cite{Mo} and we state and prove their analogue for the generalized Shioda--Inose structures of order 3. If a projective K3 surface admits an order 3 automorphism $\sigma_X$ such that the minimal resolution of $X/\sigma_X$ is isomorphic to $Km_3(A)$ for a certain Abelian surface $A$ such that $T_A\simeq  T_X$, we will say that $X$ admits a generalized Shioda--Inose structure of order 3. 
The main results of the paper are Theorems \ref{theorem: SI 1 order 3}, \ref{theorem: SI 2 order 3} and Corollaries \ref{cor: dim fami Shioda inose order 3}, \ref{cor: SI order 3 implies Si order 2 } and are summarized as follows.

\begin{theorem}\label{Main theroem intro} Let $X$ be a projective K3 surface. Then the following are equivalent:
	 \begin{itemize}\item $X$ admits a generalized Shioda--Inose structure of order 3;
		\item
    there exists a primitive embedding $T_X\hookrightarrow U\oplus A_2$;\item
    there exists a primitive embedding $\left(E_6(-1)^{\oplus 3}\right)'\hookrightarrow NS(X)$, where $\left(E_6(-1)^{\oplus 3}\right)'$ is a specific overlattice of $E_6(-1)^{\oplus 3}$, see Definition \ref{def: 3E6' and UA23E6'''};
    \item
    there exists an order 3 automorphism $\sigma_X$ on $X$ such that $\sigma_X^*$ permutes three orthogonal copies of $E_6(-1)$ contained in $\left(E_6(-1)^{\oplus 3}\right)'\hookrightarrow NS(X)$.
\end{itemize}

    The family of the K3 surfaces admitting a generalized Shioda--Inose structure of order 3 has countably many components all of dimension 1.
     
    If a K3 surface admits a generalized Shioda--Inose structure of order 3, then it admits also a Shioda--Inose structure (of order 2).\end{theorem}

The rest of the paper is devoted to the construction of specific examples: one describes certain K3 surfaces $X$ admitting an automorphism of order 3, $\sigma_X$, such that $X$ admits a generalized Shioda--Inose structure of order 3 with the automorphism $\sigma_X$. By Theorem \ref{Main theroem intro} we know that $X$ admits also an involution $\iota_X$ such that $X$ admits a Shioda--Inose structure (of order 2) with the automorphism $\iota_X$. We describe both the automorphisms $\sigma_X$ and $\iota_X$ on specific K3 surfaces. 

All the examples considered are constructed as follows: $X$ is a K3 surface endowed with an elliptic fibration $\mathcal{E}:X\ra \mathbb{P}^1$ which admits a 3-torsion section. The translation by this 3-torsion section is a symplectic automorphism of order 3 on $X$ and in all the cases considered we show that it permutes cyclically three orthogonal copies of $E_6(-1)$ embedded in $NS(X)$ and hence it is the automorphism of a generalized Shioda--Inose structure of order 3 on $X$.
 Moreover the elliptic fibration $\mathcal{E}:X\ra \mathbb{P}^1$ admits an involution of the basis, switching pairs of isomorphic fibers and preserving exactly two fibers. 
 Composing this involution with an involution acting trivially on the basis and as the elliptic involution on the fibers, one obtains a symplectic involution on $X$. Considering its action on $NS(X)$ one observes that it switches two orthogonal pairs of $E_8(-1)$ and hence it is the involution of a Shioda--Inose structure of order 2 on $X$.

Since in the examples considered the automorphism of order 3 is induced by a 3-torsion section, we are generalizing the examples constructed in \cite{vGS}, \cite{K} and \cite{Sc}. 

We discuss three different examples: the one considered in Section \ref{sec: example picard number 20} is about a very classical K3 surface $X_3$, already considered in \cite{SI} in relation of its Shioda--Inose structure of order 2 and in \cite{GP} in relation of its special automorphism of order 3;
the ones considered in Section \ref{sec: examples with Picard number 19} are families of maximal dimension (i.e. 1-dimensional) of K3 surfaces admitting a generalized Shioda--Inose structure of order 3. 
 
The main difference in the examples considered is in the relation between $T_X\simeq T_A$ and $T_{Km_3(A)}$. Due to the cohomological action of $\sigma_A^*$ on $H^2(A,\mathbb{Z})$ (which is not the identity as in the case of the involution $\iota_A^*$), this relation is not merely a multiplication of the bilinear form and in the three examples considered we indeed observe three different relations between $T_A\simeq T_X$ and $T_{Km(A)}$, as shown in the following Table
$$\begin{array}{|c|c|c|}\hline
	\mbox{Example considered in Section} &T_A\simeq T_X&T_{Km_3(A)}\\
	\hline
	\mbox{\ref{sec: example picard number 20}}&A_2&A_2\simeq T_A\simeq T_X\\
	\hline
	\mbox{\ref{subsec: 2IV*I6}}&A_2\oplus \langle -2\rangle&A_2\oplus \langle-6\rangle\\
	\hline
	\mbox{\ref{subsec: I18}}&U\oplus \langle 2\rangle&U(3)\oplus \langle 6\rangle\simeq T_A(3)\\
	\hline
\end{array}$$
 
\textbf{Acknowledgments} {\it The authors thank Bert van Geemen for many suggestions on the preliminary version of the paper and Benedetta Piroddi for useful remarks on the embeddings of some lattices.} 

\section{The main theorems}\label{sec: generalized Shioda-Inose}
\subsection{Preliminaries}	
We will apply several results on lattice theory and here we recall the main definitions used in the following.
A lattice is a free $\Z$-module of finite rank, endowed with a symmetric bilinear form $b_L:L\times L\ra \Z$. A lattice is denoted $(L,b_L)$, or simply $L$ when its bilinear form is clear by the context. 
A lattice $(L,b_L)$ is called even if for every $\ell \in L$, $b_L(\ell,\ell)\in 2\Z$. 
The dual of the lattice $L$ is $$L^*=Hom(L,\Z)\simeq \{m\in L\otimes \Q\mbox{ such that }b_L(m,\ell)\in\Z\ \forall \ell\in L\}$$
where $b_L$ denotes also the $\Q$-linear extension of $b_L$ to $L\otimes \Q$.
The discriminant group of $L$ is $A_L:=L^*/L.$ 
It is a finite group, product of cyclic groups. The minimal number of generators of $A_L$ is called length of $L$.
If $A_L$ is trivial, the lattice $L$ is said to be unimodular.
The bilinear form $b_L$ induces a bilinear form on $A_L$, called discriminant form.

We say that a lattice is $p$-elementary if $A_L\simeq (\Z/p\Z)^l$ for a certain $l$, which is the length of $L$.

Chosen a basis of $L$, the symmetric bilinear form $b_L$ is represented by a symmetric matrix, with integer entries. The signature of $L$ is the signature of this matrix, interpreted as real matrix.

The discriminant of $L$, $d(L)$, is the determinant of a matrix which represents $b_L$ with respect to a certain basis. It is independent of the choice of the basis and its absolute value  coincides with the order of the group $A_L$.

Two lattices $(L,b_L)$, $(M,b_M)$ are said to be isometric if there a isomorphism $\varphi:L\ra M$ such that $b_M(\varphi(\ell_1),\varphi(\ell_2))=b_L(\ell_1,\ell_2)$ for all $\ell_1,\ \ell_2\in L$.

A lattice $(L,b_L)$ is embedded $(M,b_M)$ if there exists an injective  homomorphism $\varphi:L\ra M$ such that $b_M(\varphi(\ell_1),\varphi(\ell_2))=b_L(\ell_1,\ell_2)$. We will say that an embedding is primitive if $M/\varphi(L)$ is torsion free. We will say that an embedding is of finite index if $M/\varphi(L)$ is a torsion group, i.e. if the ranks of $L$ and $M$ are the same.
 
In the following we use the following notation and facts:
\begin{itemize}
	\item $U$ is the unique even unimodular lattice of signature $(1,1)$;
	\item the second cohomology group of any 2-dimensional torus is isometric to the unique even unimodular lattice of signature $(3,3)$, i.e. to $U^{\oplus 3}$;
	\item the lattices $A_n$, $D_m$, $E_r$, $r=6,7,8$, are the so called ADE lattices, see \cite[Chapter I Lemma 2.12]{BHPV}. These are positive definite even lattices, spanned by their roots, i.e. by elements of self intersection 2. In particular, $A_2$ and $E_6$ are 3-elementary lattices of length 1; for any $n$, $D_{2n}$ is a  2-elementary lattice of length 2; $E_8$ is unimodular;
	\item given a lattice $L$, we denote $L(m)$ the lattice obtained by considering the $\Z$-module $L$ and by multiplying its bilinear form by $m$;
	\item $\Lambda_{K3}$ is the unique even unimodular lattice of signature $(3,19)$. The second cohomology group of any K3 surface is isometric to this lattice, and it is isometric to $U^{\oplus 3}\oplus E_8(-1)^{\oplus 2}$;
	\item $K_{12}$ is the opposite of the Coxeter--Todd lattice, defined in \cite{CT}. It is a negative definite even 3-elementary lattice of length $6$;
	\item $u(n)$ is the discriminant form of the lattice $U(n)$.
\end{itemize}

\subsection{Kummer and generalized Kummer surfaces}\label{subsec: Kummer}	
    In this section we consider special quotients of 2-dimensional tori. Given a finite order automorphism $\alpha$ on a 2-dimensional torus $A$, the Kodaira dimension of $A/\alpha$ is either $-\infty$ or $0$. If it is zero, $A/\alpha$ is birational either to a torus or to a K3 surface. We are interested in the latter case and we will call (generalized) Kummer surface a K3 surface obtained as minimal resolution of the quotient of a torus by an automorphism . We recall here the main definitions and results about this construction.

	Let $A$ be a 2-dimensional torus. Since $A$ admits a group law, the involution $\iota_A:A\ra A$ such that $\iota(a)=-a$ for each $a\in A$ is well defined.
	It has 16 fixed points, which are the 16 points in the 2-torsion group $A[2]:=\{a\in A\mbox{ such that }a+a=0_A\}$ where $+$ is the group law on $A$ and $0_A$ is its zero.
	
	As a consequence $A/\iota_A$ is a surface whose singular locus consists of 16 isolated points. Each of these singular points is a singularity of type $A_1$ and hence its minimal resolution is obtained by blowing it up once.
	
	\begin{definition}\label{defi: Kummer order 2} Let $A$ be a torus. A surface $S$ is called the Kummer surface of $A$ if it is the minimal resolution of $A/\iota_A$, where $\iota_A$ is the involution defined above. In this case we denote $S$ with $Km_2(A)$.\end{definition}
	
	We will denote $K_i$ the 16 exceptional curves arising by the resolution of the singularities of $A/\iota_A$. Their classes (still denoted by $K_i$) in the N\'eron--Severi group of $Km_2(A)$ span a sublattice of $NS(Km_2(A))$ and their intersection properties are the followings $K_i^2=-2$ and $K_iK_j=0$ if $i\neq j$. So the lattice spanned by $\{K_i\}_{i}$ is isometric to $A_1(-1)^{\oplus 16}$. This lattice is not primitive in $NS(Km_2(A))$.
	\begin{definition}\label{defi: Kummer lattice order 2} The minimal primitive sublattice contained in $NS(Km_2(A))$ and containing all the classes $K_i$ is called Kummer lattice and it is denoted by $K_{\Z/2\Z}$.
	\end{definition}    
	
	The following results on Kummer surfaces and Kummer lattice are very well known
	\begin{theorem}\label{theorem: recap Km oreder 2}{\rm(See \cite{NKummer} and \cite[Section 4]{Mo})} 
		\begin{enumerate}\item A Kummer surface $S$ is a K3 surface.\item 
The Kummer lattice $K_{\Z/2\Z}$ is an even negative definite lattice of rank 16 and discriminant form $u(2)^3$.\item 
A K3 surface $S$ is a Kummer surface if and only if the lattice $K_{\Z/2\Z}$ is primitively embedded in $NS(S)$.
The Kummer surface $Km_2(A)$ is projective if and only if the torus $A$ is projective, i.e. if $A$ is an Abelian surface.
\item A K3 surface is a Kummer surface if and only if its transcendental lattice can be primitively embedded in $U(2)^{\oplus 3}$.
\item The transcendental lattices $T_A$ and $T_{Km_2(A)}$ of a 2-dimensional torus and of its Kummer surface are related by $T_{Km_2(A)}\simeq T_A(2)$. \end{enumerate}\end{theorem}

So one obtains a projective K3 surface as quotient of an Abelian surface by an order 2 automorphism.
It is natural to ask if one can apply a similar construction for higher order automorphisms on an Abelian surface. In \cite[Section 3]{Fu} the author classified the automorphisms $\alpha$ on two dimensional tori $A$ such that the minimal resolution of $A/\alpha$ is a K3 surface. He proved that if $A$ is a torus and $\alpha\in \Aut(A)$ is such that the minimal resolution of $A/\alpha$ is a K3 surface, then the order of $\alpha$ is $|\alpha|\in\{2,3,4,6\}$, see \cite[Lemma 3.3]{Fu}. Moreover, if $|\alpha|>2$, then the torus $A$ is non generic, i.e. not all tori admit an automorphism as required with order bigger than 2.

If $A$ is a 2-dimensional torus admitting an automorphism $\sigma_A$ of order 3 such that the minimal resolution of $A/\sigma_A$ is a K3 surface, then $\sigma_A$ fixes 9 isolated points in $A$ (cf. \cite[Lemma 3.19]{Fu}, see also \cite{Be}). As a consequence, $A/\sigma_A$ has 9 singular points, all of type $A_2$. We call $J_h^{(k)}$, $k=1,\ldots, 9$, $h=1,2$ the curves arising from the desingularization of these points, in particular $J_1^{(k)}$, $J_2^{(k)}$ are mapped to the same singular point and the intersection properties of these curves are $$J_1^{(k)}J_2^{(k)}=1,\ \  (J_1^{(k)})^2=-2, \ \ J_h^{(k)}J_m^{(f)}=0,\ k\neq f.$$
The classes $J_h^{(k)}$, $k=1,\ldots, 9$, $h=1,2$ span the lattice $A_2(-1)^{\oplus 9}$ in the N\'eron--Severi group of the minimal resolution of $A/\sigma_A$ and this lattice is not primitively embedded in the N\'eron--Severi group.

The following definitions and results are analogous to the ones described above for the order 2 automorphisms $\iota_A$ on a torus $A$, (c.f. Definitions \ref{defi: Kummer order 2}, \ref{defi: Kummer lattice order 2} and Theorem \ref{theorem: recap Km oreder 2} respectively).

	\begin{definition} A K3 surface $S$ is called generalized Kummer surface of order 3 of $A$ if it exists an order three automorphism $\sigma_A$ on $A$ such that $S$ is the minimal resolution of $A/\sigma_A$. In this case we will denote $S$ as $Km_3(A)$.\end{definition}

\begin{definition} The minimal primitive sublattice contained in $NS(Km_3(A))$ and containing all the classes $J_h^{(k)}$ is called generalized Kummer lattice of order 3 and it is denoted by $K_{\Z/3\Z}$.
\end{definition}    

\begin{theorem}\label{theorem: recap Kummer order 3}\begin{enumerate}
		\item {\rm (See \cite[Section 3]{Fu}, \cite[Section 1.2]{Ba})} There exist tori $A$ which admit an order 3 automorphism $\sigma_A$ such that the minimal resolution of $A/\sigma_A$ is a K3 surface, $Km_3(A)$. In particular, $A$ admits such an automorphism if and only if the transcendental lattice $T_A$ of $A$ can be primitively embedded in $U\oplus A_2$.
		\item {\rm (See \cite[Proposition 2.3]{Be} and \cite{Ba})} The generalized Kummer lattice of order 3, $K_{\Z/3\Z}$, is an even  negative definite 3-elementary lattice of rank 18 and discriminant form which is the opposite to the one of $U(3)\oplus A_2$.
		\item {\rm (See \cite[Theorem 2.5]{Be}} A K3 surface $S$ is a generalized Kummer surface of order 3 if and only if the lattice $K_{\Z/3\Z}$ is primitively embedded in $NS(S)$.  The Kummer surface $Km_3(A)$ is projective if and only if the torus $A$ is projective, i.e. if $A$ is an Abelian surface.
	\item {\rm (See \cite[Theorem 2.5]{Be} and \cite[Section 1.3]{Ba})} A K3 surface is a generalized Kummer surface of order 3 if and only if its transcendental lattice is primitively embedded in $U(3)\oplus A_2$. \end{enumerate}
\end{theorem}

Comparing Theorem \ref{theorem: recap Km oreder 2} and Theorem \ref{theorem: recap Kummer order 3} one observes that the results are similar, but point (5) of Theorem \ref{theorem: recap Km oreder 2} states also a clear and easy relation between the transcendental lattice of $Km_2(A)$ and the one of $A$. In the case of order 3 we didn't write explicitly this relation, even if it exists, because in order to write it one needs more information. Indeed to compute the relation between $T_{Km_2(A)}$ and $T_A$ one applies the map $(\pi_2)_*:H^2(A,\Z)\ra H^2(Km_2(A),\Z)$ where $\pi_2$ is the rational $2:1$ map $A\ra Km_2(A)$. This map multiplies the form of $T_A$ by 2.
Similarly, to obtain the relation between $T_{Km_3(A)}$ and $T_A$ one has to consider the map $(\pi_3)_*$ induced by the rational $3:1$ map $A\ra Km_3(A)$. This map is studied by Barth in \cite{Ba}, but it is more complicated than $(\pi_2)_*$ and it is not just a multiplication of the form. In particular, $(\pi_3)_*(T_A)$  depends also on the embedding of $T_{A}$ in $H^2(A,\Z)$. If one determines this embedding, the relation between $T_{Km_3(A)}$ and $T_A$ is known. In Theorem \ref{theorem: SI 1 order 3} we will see how the map $(\pi_3)_*$ works and in Sections \ref{sec: example picard number 20} and \ref{sec: examples with Picard number 19} we apply it to specific examples.

\begin{rem}\label{rem: autom order 6 on A}{\rm 
	If an Abelian surface $A$ admits an order 3 automorphism $\sigma_A$ such that the minimal resolution of the quotient $A/\sigma_A$ is the K3 surface $Km_3(A)$, then it admits also an order 6 automorphism $\tau_A$ such that the minimal resolution of the quotient $A/\tau_A$ is the K3 surface $Km_6(A)$: one can chose $\tau$ to be  $\sigma_A\circ \iota_A$, see e.g. \cite[Proposition 3.5]{Fu}.}
\end{rem}

\subsection{Symplectic automorphisms on K3 surfaces}\label{subsec: aymplectic autom}
Until now we considered K3 surfaces obtained as quotients of tori by specific automorphisms. Nevertheless it is also possible to obtain K3 surfaces as quotient of other K3 surfaces by specific automorphisms. In this section we discuss this situation recalling the main results about this construction.

\begin{definition}
	An automorphism $\alpha$ on a K3 surface $S$ is symplectic if it acts trivially on the two holomorphic forms on $S$,  i.e. if $\alpha^*_{|H^{2,0}(S)}=\id_{|H^{2,0}(S)}$.
	\end{definition}	
The finite order symplectic automorphisms on K3 surfaces are studied intensively and we recap here the  main results. 
\begin{theorem}\label{theorem: recap 1 sympl autom}{\rm (See \cite{NSympl})} Let $\alpha$ be an order $n$ symplectic automorphism on a K3 surface $S$ and $Y$ be the minimal resolution of $S/\alpha$. Then:
\begin{enumerate}
	\item $2\leq n\leq 8$;
	\item $Y$ is a K3 surface;
	\item if $n=2$, then $X/\alpha$ has 8 singularities of type $A_1$;
	\item if $n=3$, then $X/\alpha$ has 6 singularities of type $A_2$.
\end{enumerate}	
\end{theorem}
In \cite{NSympl}, it is also proved that there are two lattices which characterize the K3 surfaces which admit a symplectic automorphism and the ones which are obtained as minimal resolution of the quotient of a K3 surface by an automorphism. Let us define these lattices.

If $S$ is a K3 surface and $\alpha_n$ is an order $n$ symplectic automorphism on $S$, $S/\alpha_n$ is singular in some isolated points $p_k$. Denoted $Y_n$ the minimal resolution of $S/\alpha_n$, for each of the singular point $p_k\in S/\alpha_n$ there are some rational curves $M_h^{(k)}\subset Y_n$ arising from the desingularization of the point $p_k$. The number and the configuration of these curves depend on $n$, see e.g. Theorem \ref{theorem: recap 1 sympl autom} for the cases $n=2,3$. The curves $M_h^{(k)}$ correspond to the classes in $NS(Y_n)$, which span a sublattice of the N\'eron--Severi lattice. This lattice is not primitive in $NS(Y_n)$.
\begin{definition}\label{defi: Mn}
 We denote $M_{\Z/n\Z}$ the minimal primitive sublattice of $NS(Y_n)$ containing all the classes of the curves $M_h^{(k)}$. In particular, if $n=2$, the lattice $M_{\Z/2\Z}$ is called Nikulin lattice and it is often denoted by $N$. 
\end{definition}
	
	In \cite{NSympl} Nikulin also proved  that the action of $\alpha_n^*$ on $H^2(S,\mathbb{Z})$ is essentially unique and in particular the isometry class of the lattices $(H^2(S,\mathbb{Z}))^{\alpha^*_n}$ and $\left((H^2(S,\mathbb{Z}))^{\alpha^*_n}\right)^{\perp}$ depends only on $n$. Hence the following definition is well posed since $H^2(S,\mathbb{Z})$ is isometric to the standard lattice $\Lambda_{K3}$ which does not depend on $S$ and $\alpha_n^*$ is the isometry induced on $\Lambda_{K3}$ by the isometry $\Lambda_{K3}\simeq H^2(S,\Z)$.

\begin{definition}\label{defi: Omegan}
	We denote $\Omega_n$ the lattice $\left(\Lambda_{K3}^{\alpha^*_n}\right)^{\perp}$.
\end{definition}	

\begin{theorem}\label{theorem: recap 2 sympl K3} 
	\begin{enumerate}
		\item {\rm (See \cite[Lemma 4.2]{NSympl})} A K3 surface $S_n$ admits a symplectic automorphism of order $n$ if and only if $\Omega_n$ is primitively embedded in $NS(S_n)$.
		\item {\rm (See \cite{NSympl}, \cite[Theorem 5.2]{G2})} A K3 surface $Y_n$ is the minimal resolution of the quotient of a K3 surface by an order $n$ symplectic automorphism if and only if $M_{\Z/n\Z}$ is primitively embedded in $NS(Y_n)$.
		\item {\rm (See \cite[Proposition 7.1 and Lemma 10.2]{NSympl})} For each $2\leq n\leq 8$, the lattices $M_{\Z/n\Z}$ and $\Omega_n$ are negative definite lattices such that $\rk(M_{\Z/n\Z})=\rk(\Omega_n)$.
		\item {\rm (See \cite[Theorem 5.7]{Mo})} If $n=2$, then $\rk(\Omega_2)=\rk(M_{\Z/2\Z})=8$; $\Omega_2$ is isometric to $E_8(-2)$, it is a 2-elementary lattice and its discriminant form is $u(2)^4$; $M_{\Z/2\Z}$ is a 2-elementary lattice and its discriminant form is $u(2)^3$.
		\item {\rm (See \cite[Proposition 7.1]{NSympl}, \cite[Theorem 4.3]{GSprime})} If $n=3$, then $\rk(\Omega_3)=\rk(M_{\Z/3\Z})=12$; $\Omega_3$ is isometric to $K_{12}$, it is a 3-elementary lattice and its discriminant form is the one of $U(-3)^{\oplus 2}\oplus A_2^{\oplus 2}$; $M_{\Z/3\Z}$ is a 3-elementary lattice and its discriminant form is the one of $U(-3)\oplus A_2^{\oplus 2}$.
		
	\end{enumerate}	
\end{theorem}

Due to \cite[Theorem 1.14.2]{Nforms} the embeddings of the lattices $M_{\Z/n\Z}$ and $\Omega_n$, $n=2,3$ in $\Lambda_{K3}$ are unique up to isometries of $\Lambda_{K3}$, hence in the following we are free to fix specific embeddings.  

As mentioned above, the action induced on the second cohomology group by a symplectic automorphism on a K3 surface is essentially unique, \cite[Theorem 4.7]{NSympl}.

In \cite{Mo}, Morrison describes explicitly this isometry in the case $n=2$.
\begin{theorem}\label{theorem: action order three}{\rm (See \cite{Mo})} Let $S$ be a K3 surface admitting a symplectic involution $\iota$. Then there exists an isometry $\varphi: H^2(S,\mathbb{Z})\rightarrow U^{\oplus 3}\oplus E_8(-1)^{\oplus 2}$ such that $\varphi\circ\iota^*\circ \varphi^{-1}$ acts as the identity on $U^{\oplus 3}$ and switches the two copies of $E_8(-1)$. 
\end{theorem}

Similarly in \cite{GP} the explicit description of the isometry is provided in case $n=3$.

In order to state the result for $n=3$, we need some definitions.
\begin{definition}\label{def: 3E6' and UA23E6'''}
	Let $e_i$ be a basis of $E_6(-1)$ with the following intersection properties:\begin{eqnarray}\label{eq:E6}\xymatrix{e_1\ar@{-}[r]&e_2\ar@{-}[r]&e_3\ar@{-}[r]\ar@{-}[d]&e_4\ar@{-}[r]&e_5\\&&e_6}\end{eqnarray}

	Let $v=(e_1+2e_2+e_4+2e_5)/3$.
	
 The lattice $(E_6(-1)^{\oplus 3})'$ is the overlattice of index 3 of $E_6(-1)^{\oplus 3}$ obtained adding the vector $x=v^{(1)}+v^{(2)}+v^{(3)}$ where $e_i^{(j)}$, $i=1,\ldots, 6$ is the basis of the $j$-th copy of $E_6(-1)$ in $E_6(-1)^{\oplus 3}$ and $v^{(j)}$ is the corresponding vector $v$ as above.
 
 We denote $(U\oplus A_2\oplus E_6(-1)^{\oplus 3})''$ the overlattice of index $3^2$ obtained by adding to $U\oplus A_2\oplus E_6(-1)^{\oplus 3}$ the vectors $x$ and $y=(a_1+2a_2)/3+ v^{(1)}-v^{(2)}$ where $a_i$, $i=1,2$ is a basis of $A_2$ such that $a_i^2=2$, $a_1a_2=-1$. 
 
 Equivalently  $(U\oplus A_2\oplus E_6(-1)^{\oplus 3})''$ is the overlattice of index 3 obtained by adding to $U\oplus A_2\oplus (E_6(-1)^{\oplus 3})'$ the vector $y$.
\end{definition}

The lattice $(U\oplus A_2\oplus E_6(-1)^{\oplus 3})''$ is an even unimodular lattice of signature $(3,19)$. This suffices to conclude that $(U\oplus A_2\oplus E_6(-1)^{\oplus 3})''$ is isometric to $\Lambda_{K3}$. Due to \cite[Theorem 1.14.2]{Nforms} the embedding of $(E_6(-1)^{\oplus 3})'$ in $\Lambda_{K3}$ is unique up to isometries.

\begin{theorem}{\rm (See \cite[Section 3.1]{GP})}
	Let $S$ be a K3 surface admitting an order 3 symplectic automorphism $\sigma$. Then there exists an isometry $$\varphi: H^2(S,\mathbb{Z})\rightarrow (U\oplus A_2\oplus E_6(-1)^{\oplus 3})''$$ such that $\varphi\circ\iota^*\circ \varphi^{-1}$ acts as the identity on $U\oplus A_2$ and permutes cyclically the three copies of $E_6(-1)$. 
\end{theorem}

\subsection{Shioda--Inose and generalized Shioda--Inose structures}\label{subsec: Shioda Inose}
One of the aims of the paper \cite{Mo} was to relate K3 surfaces and Abelian surfaces which admit the same K3 surface as quotient by an order 2 automorphism. To do that, the author introduced the so-called Shioda--Inose structures. Here we recall his results and then we generalize them to the order 3 case.
\subsubsection{The definitions}
\begin{definition}\label{def: SI order 2}
	A Shioda Inose structure (of order 2) is a triple $(A,X,\iota_X)$ such that:
	\begin{enumerate}\item $A$ is a 2-dimensional torus;
	\item $X$ is a K3 surface and $\iota_X$ is a symplectic involution on $X$;
\item the minimal resolution of $X/\iota_X$ is isomorphic to $Km_2(A)$;
\item the transcendental lattices of $A$ and $X$ are isometric, i.e. $T_A\simeq T_X$.\end{enumerate}
\end{definition}

\begin{definition}\label{def: SI order 3} A generalized Shioda--Inose structure of order 3 is a quadruple $(A,\sigma_A, X,\sigma_X)$ where:
	\begin{enumerate}\item $A$ is a 2-dimensional torus admitting a symplectic automorphism $\sigma_A$ of order 3 such that the minimal resolution of $A/\langle\sigma_A\rangle$ is a K3 surface $Km_3(A)$; \item $X$ is a K3 surface and $\sigma_X$ is an order 3 symplectic automorphism on $X$;
		\item the minimal resolution of $X/\langle\sigma_X\rangle$ is isomorphic to $Km_3(A)$; \item the transcendental lattices of $A$ and $X$ are isometric, i.e. $T_A\simeq T_X$.\end{enumerate}\end{definition}

\begin{rem}\label{rem: generalization definition of SI order n}{\rm One could generalize the previous definitions of Shioda--Inose structures to the order $n$ cases, for $n=2,3,4,6$. In this case one could define a Shioda--Inose structure of order $n$ to be the the quadruple $(A,\tau_A,X,\tau_X)$ where: $A$ is an Abelian surface and $\tau_A\in Aut(A)$; $X$ is a K3 surface and $\tau_X\in Aut(X)$; $|\tau_A|=|\tau_X|=n$; the minimal resolutions of $A/\tau_A$ and of $X/\tau_X$ are isomorphic to $Km_n(A)$; $T_X\simeq T_A$.

Nevertheless we have no results which guarantee the existence of these structures for $n\neq 2,3$. We have an explicit example in which the generalized Shioda--Inose structure for order 6 can not be constructed, see Section \ref{subsubsection: remarks n=6}, but we have not a characterization of the Abelian surfaces for which is possible/impossible to construct a generalized Shioda--Inose structure as defined above.}\end{rem}

With the previous notation a Shioda-Inose structure of order $n$ is represented by the following diagram
$$\xymatrix{A\ar[d]\ar@{-->}[dr]
	&&X\ar[d]\ar@{-->}[dl]\\
	A/\tau_A&Km_n(A)\ar[r]\ar[l]&X/\tau_X}$$
where the horizontal arrows are birational morphisms, the vertical ones are quotient maps of order $n$ and hence the diagonal ones are generically $n:1$.
Moreover one requires $T_A\simeq T_X$. If $|\tau|=2$ we often denote it with $\iota$, if $|\tau|=3$ with $\sigma$.

\subsubsection{The case of order 2: recap of known results}
\begin{theorem}\label{theorem: recap 1 SI order 2}{\rm (See \cite[Theorem 5.7]{Mo})} Let $X$ be a K3 surface such that $E_8(-1)^{\oplus 2}$ is primitively embedded in $NS(X)$. Then there exists a symplectic involution $\iota_X$, such that denoted by $Y$ the minimal resolution of $X/\iota_X$ and by $\pi:X\dashrightarrow Y$ the rational quotient map, the followings hold:
	\begin{enumerate}
		\item the lattice $M_{\Z/2\Z}\oplus E_8(-1)$ is primitively embedded in $NS(Y)$;
		\item $\pi_*$ induces a Hodge isometry $T_X(2)\simeq T_Y$;
		\item $M_{\Z/2\Z}\oplus E_8(-1)$ is a negative definite rank 16 lattice whose discriminant form is $u(2)^3$. 
	\end{enumerate}
\end{theorem}
\begin{rem}\label{rem: same discriminant form N+E8 and K}{\rm 
By Theorems \ref{theorem: recap Km oreder 2} and \ref{theorem: recap 1 SI order 2} it follows that both $K_{\Z/2\Z}$ and $M_{\Z/2\Z}\oplus E_8(-1)$
are rank 16 negative definite 2-elementary lattices whose discriminant forms are $u(2)^3$. Nevertheless they are not isometric, cf. Proposition \ref{prop: negative definite lattices and genus}, see also \cite{Mo}. We will generalize this result for the order 3 case.}\end{rem}

\begin{theorem}\label{theorem: recap 2 SI order 2}{\rm (See \cite[Theorem 6.3]{Mo})} Let $X$ be a projective K3 surface. Then the followings are equivalent:
	\begin{enumerate}
		\item there exist an Abelian surface $A$ and a symplectic involution $\iota_X\in Aut(X)$ such that $(A,X,\iota_X)$ is a Shioda--Inose structure;
		\item there exists an Abelian surface $A$ such that $T_A\simeq T_X$;
		\item there exists a primitive embedding $T_X\hookrightarrow U^{\oplus 3}$;
		\item there exists a primitive embedding of $E_8(-1)^{\oplus 2}\hookrightarrow NS(X)$.
	\end{enumerate}
\end{theorem}

\subsubsection{The case of order 3: new results} 
In order to generalize Morrison's results to the order 3 case, one needs some preliminary results and definitions.

We defined the lattices $M_{\Z/n\Z}$ in Definition \ref{defi: Mn}, now we give an explicit basis of this lattice if $n=3$, see \cite[Section 6]{NSympl}.
Let $M_1$ and $M_2$ be a basis of $A_2(-1)$ such that $M_1M_2=1$, $M_i^2=-2$, and $z:=(M_1+2M_2)/3\in A_{A_2(-1)}$. 
In the lattice $A_2(-1)^{\oplus 6}$ we denote $M_1^{(j)}$, $M_2^{(j)}$, $j=1,\ldots, 6$ the basis of the $j$-th copy of $A_2(-1)$ and by $z^{(j)}:=(M_1^{(j)}+2M_2^{(j)})/3$.
The lattice $M_{\Z/3\Z}$ is the overlattice of index 3 of $A_2(-1)^{\oplus 6}$ generated by the classes $M_i^{(j)}$, $i=1,2$, $j=1,\ldots, 6$ and by the class $\hat{M}:=\sum_{j=1}^{6}z^{(j)}.$

\begin{definition}\label{def: overlatticeM3E6}
We denote $(M_{\Z/3\Z}\oplus E_6(-1))'$ the overlattice of index 3 of $M_{\Z/3\Z}\oplus E_6(-1)$ obtained by adding the vector $n:=(e_1+2e_2+e_4+2e_5)/3-z^{(1)}+z^{(3)}-z^{(4)}+z^{(5)}$
where $e_i$, $i=1,\ldots, 6$ is the basis of $E_6(-1)$ defined in \eqref{eq:E6}.
\end{definition}

The lattice $(M_{\Z/3\Z}\oplus E_6(-1))'$ is primitively embedded in $\Lambda_{K3}$ (see \cite[Section 3.4]{GP}) and this embedding is unique up to isometry by \cite[Theorem 1.14.2]{Nforms}.

We consider the basis $\{u_1,u_2,a_1,a_2\}$ of the lattice $U\oplus A_2$, so that $u_i^2=0$, $a_i^2=2$, $i=1,2$, $u_1u_2=a_1a_2=-1$, the other intersections are zero.

We consider the basis $\{w_1,w_2,b_1,b_2\}$ of the lattice $U(3)\oplus A_2(3)$, so that $w_i^2=0$, $b_i^2=6$, $w_1w_2=b_1b_2=-3$, the other intersections are zero. 

We observe that $\{(b_1+2b_2)/3, b_2\}$ is a basis of $A_2$ (indeed the intersection form computed on $(b_1+2b_2)/3$ and $(b_1+2b_2)/3-b_2$ gives the standard matrix of $A_2$).
\begin{definition}\label{def: gamma}
	Let
	$$\gamma:U\oplus A_2\ra U(3)\oplus A_2\subset \left(U(3)\oplus A_2(3)\right)\otimes \Q$$ be the map such that $\gamma(u_i)=w_i\ \ \gamma(a_1)=(b_1+2b_2)/3,\ \gamma(a_2)=b_2.$
\end{definition}

The following is the order 3 analogue of Theorem \ref{theorem: recap 1 SI order 2}.

\begin{theorem}\label{theorem: SI 1 order 3} Let $X$ be a K3 surface such that $\left(E_6(-1)^{\oplus 3}\right)'$ is primitively embedded in $NS(X)$. Then there exists an order 3 symplectic automorphism $\sigma_X$ such that, denoted by $Y$ the minimal resolution of $X/\sigma_X$ and by $\pi:X\dashrightarrow Y$ the rational quotient map, the following hold:
	\begin{enumerate}
		\item the lattice $\left(M_{\Z/3\Z}\oplus E_6(-1)\right)'$ is primitively embedded in $NS(Y)$;
		\item the transcendental lattice of $X$ is primitively embedded in $U\oplus A_2$ and the map $\pi_*$ acts on $T_X\subset U\oplus A_2$ as the restriction of $\gamma$ to $T_X$.
		\item $\left(M_{\Z/3\Z}\oplus E_6(-1)\right)'$ is a rank 18 negative definite 3-elementary lattice whose discriminant form is the one of $U(-3)\oplus A_2(-1)$. 
	\end{enumerate}
\end{theorem}
\proof We define an isometry $r$ of $H^2(X,\Z)$ which has order three: it permutes cyclically the three copies of $E_6(-1)$ in $\left(E_6(-1)^{\oplus 3}\right)'$ and acts as the identity on the orthogonal of $\left(E_6(-1)^{\oplus 3}\right)'$ in $H^2(X,\Z)$. The lattice $(H^2(X,\Z)^r)^{\perp}$ is $K_{12}$ (see \cite[Proposition 3.1]{GP}) and it is primitively embedded in $\left(E_6(-1)^{\oplus 3}\right)'$. Since $\left(E_6(-1)^{\oplus 3}\right)'$ is primitively embedded in $NS(X)$, $K_{12}$ is primitively embedded in $NS(X)$ and hence $X$ admits an order 3 symplectic automorphism $\sigma_X$ which acts trivially on the orthogonal of $K_{12}\subset \left(E_6(-1)^{\oplus 3}\right)'$ and acts as $r$ on $K_{12}$. The isometry $\sigma_X^*$ on $\Lambda_{K3}\simeq H^2(X,\Z)$ coincides with $r$. We proved that $X$ admits an order 3 symplectic automorphism $\sigma_X$ acting as an order 3 cyclic permutation of the three copies of $E_6(-1)$ in $\left(E_6(-1)^{\oplus 3}\right)'$. Now we prove that the minimal resolution $Y$ of $X/\sigma_X$ satisfies the listed properties.

In \cite[Section 3.4]{GP} it is observed that $H^2(Y,\Z)$ is an overlattice of finite index of $M_{\Z/3\Z}\oplus \pi_*(H^2(X,\Z))$. Moreover in \cite[Section 3.5]{GP} the gluing vectors needed to obtain $H^2(Y,\Z)$ as overlattice of  $M_{\Z/3\Z}\oplus \pi_*(H^2(X,\Z))$ were determined.
Since $M_{\Z/3\Z}$ is by definition contained in $NS(Y)$, $NS(Y)$ is an overlattice of finite index of $M_{\Z/3\Z}\oplus \pi_*(NS(X))$ and the gluing vectors which allows to obtain $NS(Y)$ by $M_{\Z/3\Z}\oplus \pi_*(NS(X))$ are the among the ones considered to obtain $H^2(Y,\Z)$ by $M_{\Z/3\Z}\oplus \pi_*(H^2(X,\Z))$.

In \cite[Proposition 3.2]{GP} the map $\pi_*:H^2(X,\Z)\ra \pi_*(H^2(X,\Z))\subset H^2(Y,\Z)$ is described as extension to the overlattice $(U\oplus A_2\oplus E_6(-1)^{\oplus 3})''$ of the map:
$$\begin{array}{ccccccccccccccc}\pi_*:&U\oplus &A_2\oplus &E_6(-1)\oplus &E_6(-1)\oplus& E_6(-1)&\ra& U(3)\oplus &A_2(3)\oplus &E_6(-1)\\
	&(u,&a,&e,&f,&g)&\mapsto& (u,&a,& e+f+g).\end{array}$$
In particular, $\pi_*(E_6(-1)^{\oplus 3})\simeq E_6(-1)$ and in the proof of \cite[Proposition 3.2]{GP}  it is shown that  $\pi_*((E_6(-1)^{\oplus 3})')\simeq \pi_*(E_6(-1)^{\oplus 3})\simeq E_6(-1)$. Since $(E_6(-1)^{\oplus 3})'$ is primitively embedded in $NS(X)$, $E_6(-1)\simeq \pi_*((E_6(-1)^{\oplus 3})')$ is primitively embedded in $NS(Y)$. So $NS(Y)$ contains a finite index overlattice of $M_{\Z/3\Z}\oplus E_6(-1)$. In \cite[Proposition 3.4]{GP} it is proved that the class $n$, see Definition \ref{def: overlatticeM3E6}, is contained in $H^2(Y,\Z)$. Since it is also contained in $A_{M_{\Z/3\Z}\oplus E_6(-1)}$, the minimal primitive sublattice of $H^2(Y,\Z)$ containing $M_{\Z/3\Z}\oplus E_6(-1)$ contains also $n$. It follows that $NS(Y)$ contains $M_{\Z/3\Z}\oplus E_6(-1)$ and the vector $n$, i.e. it contains the lattice $(M_{\Z/3\Z}\oplus E_6(-1))'$. By the description of $H^2(Y,\Z)$ given in \cite[Proposition 3.4]{GP} it follows that $(M_{\Z/3\Z}\oplus E_6(-1))'$ is primitively embedded in $H^2(Y,\Z)$ and so $(M_{\Z/3\Z}\oplus E_6(-1))'$ is primitively embedded in $NS(Y)$.

Since $$\left(E_6(-1)^{\oplus 3}\right)'\hookrightarrow NS(X)\mbox{ one has }\left(\left(E_6(-1)^{\oplus 3}\right)'\right)^{\perp}\hookleftarrow (NS(X))^{\perp}\simeq T_X$$ where the perpendicular are considered in $H^2(X,\Z)$. We described $H^2(X,\Z)$ as overlattice of index 3 of $U\oplus A_2\oplus (E_6(-1)^{\oplus 3})'$ and we recall that the embedding of $(E_6(-1)^{\oplus 3})'$ in $\Lambda_{K3}\simeq H^2(X,\Z)$ is unique (up to isometries). Hence one sees that $\left(\left(E_6(-1)^{\oplus 3}\right)'\right)^{\perp}\simeq U\oplus A_2$ and one deduces that $T_X\hookrightarrow U\oplus A_2$.

The discriminant group and form of the lattice $(M_{\Z/3\Z}\oplus E_6(-1))'$ can be explicitly computed, since we have a basis of this lattice. Alternatively one can observe that the orthogonal to  $(M_{\Z/3\Z}\oplus E_6(-1))'$ in $H^2(Y,\Z)$ is $U(3)\oplus A_2$, and since $H^2(Y,\Z)$ is a unimodular lattice by \cite[1.6.1]{Nforms} it follows that the discriminant form of $(M_{\Z/3\Z}\oplus E_6(-1))'$ is the opposite of the one of $U(3)\oplus A_2$. To find the orthogonal of $(M_{\Z/3\Z}\oplus E_6(-1))'$ in $H^2(Y,\Z)$ it suffices to consider the description of $H^2(Y,\Z)$ as overlattice of $M_{\Z/3\Z}\oplus \pi_*(H^2(X,\Z)\simeq M_{\Z/3\Z}\oplus (E_6(-1))\oplus U(3)\oplus A_2$ given in \cite[Proposition 3.4]{GP}.
\endproof

We now introduce a lattice theoretic result, which will be useful in the following. Part of this result was already proved in \cite{Mo}, in particular the statement about lattice which are related with Kummer and Shioda--Inose structures of order 2, i.e. on the pair $\left(K_{\Z/2\Z},E_8\oplus M_{\Z/2\Z}\right)$, see Remark \ref{rem: same discriminant form N+E8 and K}.

\begin{prop}\label{prop: negative definite lattices and genus}
	Let $(L_1,L_2)$ be one of the following pairs of lattices: $$\begin{array}{c}\left(K_{\Z/2\Z},E_8(-1)\oplus M_{\Z/2\Z}\right),\ \ \left(K_{\Z/3\Z},(E_6(-1)\oplus M_{\Z/3\Z})'\right),\\ \left((E_6(-1)\oplus E_6(-1)\oplus E_6(-1))',E_8(-1)\oplus E_8(-1)\oplus A_2(-1)\right).\end{array}$$ 
	
	Then $L_1$ and $L_2$ have the same rank, discriminant group and form, but they are not isometric.
	
	Let $R_i$ be an even hyperbolic lattice which contains primitively $L_i$ and such that $\rk(R_i)=\rk(L_i)+1$. If $R_1$ and $R_2$ have the same discriminant group and form, then they are isometric. 
\end{prop} 
\proof The discriminant forms and lattices of $K_{\Z/2\Z}$, $M_{\Z/2\Z}$, $K_{\Z/3\Z}$ and $(E_6(-1)\oplus M_{\Z/3\Z})'$ are given in Theorems \ref{theorem: recap Km oreder 2}, \ref{theorem: recap 2 sympl K3}, \ref{theorem: recap Kummer order 3}, \ref{theorem: SI 1 order 3} respectively. Since $E_8(-1)$ is unimodular, the discriminant group and form of $M_{\Z/2\Z}$ (resp. $A_2(-1)$) are the same as the ones of $M_{\Z/2\Z}\oplus E_8(-1)$ (resp. $E_8(-1)\oplus E_8(-1)\oplus A_2(-1)$). The orthogonal complement of $(E_6(-1)\oplus E_6(-1)\oplus E_6(-1))'$ in $\Lambda_{K3}$ is $U\oplus A_2$ (see Definition \ref{def: 3E6' and UA23E6'''} and Theorem \ref{theorem: action order three}), so the discriminant form of  $(E_6(-1)\oplus E_6(-1)\oplus E_6(-1))'$ is the same as the one of $A_2(-1)$. 

In order to prove that the lattices $L_1$ and $L_2$ in the statement are not isometric, it suffices to compute the number of vectors of length $-2$ in each of them. All these lattices are negative definite, so the number of vectors of a given length is finite and a computer calculation allows to find it. In particular, one finds that this number is as in the table:
$$\begin{array}{|c|c||c|c||c|c|}\hline K_{\Z/2\Z}&E_8(-1)\oplus M_{\Z/2\Z}&K_{\Z/3\Z}&\left(E_6(-1)\oplus M_{\Z/3\Z}\right)'&\left(E_6(-1)^{\oplus 3}\right)'&(E_8^{\oplus 2}\oplus A_2)(-1)\\ \hline
32& 256&54&108& 216& 486\\ \hline\end{array}$$

The lattices $R_i$ are overlattices of finite index (possibly 1) of $\langle 2d\rangle\oplus L_i$ where $d$ is a positive integer. It suffices to check that lattices $R_i$ are unique in their genus to conclude that $R_1\simeq R_2$. By \cite[Corollary 1.13.3]{Nforms}, the signature $(s_+,s_-)$, the discriminant group and the discriminant form determines uniquely a lattice $R_i$ if both $s_+$ and $s_-$ are positive and $l(R_i)\leq \rk(R_i)-2$, where $l(R_i)$ is the length of $R_i$. Since $R_i$ is hyperbolic the condition on the signature is surely satisfied. Moreover, in all the cases considered  $l(L_i)\leq \rk(L_i)-2$ and $l(R_i)\leq l(L_i)+1$. So we have $$l(R_i)\leq l(L_i)+1\leq \rk(L_i)+1-2=\rk(R_i)-2.$$
This implies that $R_1\simeq R_2$ in all the cases considered. \endproof

In view of the second part of the statement of the previous proposition, we are forced to consider hyperbolic lattices to conclude the existence of certain isometries. Since we will be interested in the N\'eron--Severi group of surfaces, we observe that if the N\'eron--Severi group of a surface is hyperbolic then the surface is projective.
The following is the order 3 analogue of Theorem \ref{theorem: recap 2 SI order 2}.

\begin{theorem}\label{theorem: SI 2 order 3} Let $X$ be a projective K3 surface. Then the followings are equivalent:
	\begin{enumerate}
		\item there exist an Abelian surface $A$, an automorphism $\sigma_A\in \Aut(A)$ and a symplectic automorphism $\sigma_X\in Aut(X)$ such that $(A,\sigma_A,X,\sigma_X)$ is a generalized Shioda--Inose structure of order 3;
		\item there exists an Abelian surface $A$ admitting an order 3 automorphism $\sigma_A$ such that $A/\sigma_A$ is birational to a K3 surface and $T_A\simeq T_X$;
		\item there exists a primitive embedding $T_X\hookrightarrow U\oplus A_2$;
		\item there exists a primitive embedding of $\left(E_6(-1)^{\oplus 3}\right)'$ in $NS(X)$.
	\end{enumerate}
\end{theorem}
\proof $(1)\Rightarrow (2)$ follows from the definition of generalized Shioda--Inose structure.

$(2)\Rightarrow (3)$ The assumption $(2)$ implies the existence of an Abelian surface $A$ such that $T_A\simeq T_X$. Moreover $A$ admits an order 3 automorphism $\sigma_A$ such that $A/\sigma_A$ is birational to a K3 surface. By Theorem \ref{theorem: recap Kummer order 3}, it follows that $T_A$ is primitively embedded $U\oplus A_2$, so $T_X\simeq T_A$ is primitively embedded in $U\oplus A_2$.

$(3)\Rightarrow (4)$ Since $X$ is projective and $T_X\hookrightarrow U\oplus A_2$, $\rk(T_X)\leq 3$ and $sgn(T_X)=(2,t)$ with $t\leq 1$. By \cite[Theorem 1.14.4]{Nforms}, $T_X$ admits a unique primitive embedding in $\Lambda_{K3}\simeq (U\oplus A_2\oplus E_6(-1)^{\oplus 3})''$. We can assume it is embedded in the first two summands and we obtain that $NS(X)=T_X^{\perp}$ contains $(U\oplus A_2)^{\perp}\simeq(E_6(-1)^{\oplus 3})'$.

$(4)\Rightarrow (1)$ In the proof of Theorem \ref{theorem: SI 1 order 3} we have shown that if $\left(E_6(-1)^{\oplus 3}\right)'$ is primitively embedded in $NS(X)$, then $X$ admits an order 3 symplectic automorphism $\sigma_X$ such that $\sigma_X^*$ permutes cyclically the three copies of $E_6(-1)$ in $\left(E_6(-1)^{\oplus 3}\right)'$. By Theorem \ref{theorem: SI 1 order 3}, $T_X$ is primitively embedded in $U\oplus A_2$ and the embedding is uniquely determined. Denoted by $Y$ the minimal resolution of $X/\sigma_X$, its transcendental lattice $T_Y$ is isometric to $\gamma(T_X)$ (see Definition \ref{def: gamma} and Theorem \ref{theorem: SI 1 order 3}) and $(M_{\Z/3\Z}\oplus E_6(-1))'$ is primitively embedded in $NS(Y)$. By \cite[Theorem 1.4.14]{Nforms}, the lattice $NS(Y)$ is uniquely determined by its signature and discriminant form since it is hyperbolic, its rank is bigger than or equal to 19 and its length is less than or equal to 3. Since the negative definite lattices $K_{\Z/3\Z}$ and $(M_{\Z/3\Z}\oplus E_6(-1))'$ have the same discriminant form and rank, by \cite[Lemma 2.3]{Mo}, $K_{\Z/3\Z}$ is primitively embedded in $NS(Y)$. By Theorem \ref{theorem: recap Kummer order 3}, there exists an Abelian surface and an order 3 automorphism $\sigma_A$ such that $Y$ is isomorphic to the minimal resolution of $A/\sigma_A$, i.e. $Y\simeq Km_3(A)$. It remains to determine the transcendental lattice of $A$.

The transcendental lattice of $Km_3(A)$ is $(\pi_A)_*(T_A)\subset H^2(Km_3(A),\Z)$ where $\pi_A$ is the rational quotient map $A\dashrightarrow Km_3(A)$. In \cite{Ba} the action of $(\pi_A)_*$ is described, in particular it acts on $U\oplus A_2\subset U^{\oplus 3}$ as the map $\gamma$, see Definition \ref{def: gamma} and \cite[Section 1.3]{Ba}. Hence $T_{Km_3(A)}=\gamma(T_A)$ and since $Y\simeq Km_3(A)$, we have $T_Y\simeq T_{Km_3(A)}\simeq \gamma(T_A)$. We already observed that $T_Y\simeq \gamma(T_X)$, hence  $$\gamma(T_X)\simeq T_Y\simeq T_{Km_3(A)}\simeq \gamma(T_A).$$
So $T_A$ has to be embedded in $U\oplus A_2$ as $T_X$, because this guarantees that $\gamma(T_A)=\gamma(T_X)$.

Therefore we constructed an order three automorphism $\sigma_X$ on $X$; an Abelian surface $A$ with $T_A\simeq T_X$; an order 3 automorphism $\sigma_A$ on $A$ such that the minimal resolution of $X/\sigma_X$ and $A/\sigma_A$ are isomorphic. It follows that $(A,\sigma_A, X,\sigma_X)$ is a generalized Shioda--Inose structure of order 3.\endproof 

Let $\mathcal{L}_{2d}:=\{(\langle 2d\rangle\oplus (E_6(-1)^{\oplus 3})')\mbox{-polarized K3 surfaces}\}/\sim$ where $\sim$ is the isomorphisms of polarized K3 surfaces. Let $\left(\langle 2d\rangle\oplus (E_6(-1)^{\oplus 3})'\right)'$ be the unique overlattice of index 3 of $\langle 2d\rangle\oplus (E_6(-1)^{\oplus 3})'$ in which both $\langle 2d\rangle$ and $(E_6(-1)^{\oplus 3})'$ are primitively embedded. We will prove that it exists only for certain values of $d$ and that in this case it is unique.
Let $\mathcal{M}_{2d}:=\{(\left(\langle 2d\rangle\oplus (E_6(-1)^{\oplus 3})'\right)')\mbox{-polarized K3 surfaces}\}/\sim$.
\begin{corollary}\label{cor: dim fami Shioda inose order 3}
The family of the projective K3 surfaces $X$ such that there exists a generalized Shioda--Inose structure of order 3 $(A,\sigma_A,X,\sigma_X)$ is 
$$\left(\bigcup_{d\in\mathbb{N}_{>0}}\mathcal{L}_{2d}\right)\bigcup \left(\bigcup_{\footnotesize{\begin{array}{l}e\in\mathbb{N}_{>0},\ e\equiv_93\\\mbox{ or }e\equiv_96\end{array}}}\mathcal{M}_{2e}\right).$$
In particular it has countably many irreducible components of dimension 1.
\end{corollary} 
\proof By the previous Theorem, there exists a generalized Shioda--Inose structure with $X$ as K3 surface if and only if $(E_6(-1)^{\oplus 3})'$ is primitively embedded in $NS(X)$. Since $X$ is projective, there exists a class $L$ with positive self intersection in the orthogonal complement of $(E_6(-1)^{\oplus 3})'$ in $NS(X)$. Hence $NS(X)$ has rank at least 19 and if it has rank 19, then it is an overlattice of finite index (possibly 1) of $\langle 2d\rangle\oplus (E_6(-1)^{\oplus 3})'$. The discriminant form of $(E_6(-1)^{\oplus 3})'$ is $\Z_3(\frac{-2}{3})$. Arguing as in \cite[Proposition 4.2]{GP} one obtains that if $d\not\equiv 0\mod 3$, there are no proper overlattices of $\langle 2d\rangle\oplus (E_6(-1)^{\oplus 3})'$ of finite index. If $d\equiv 3\mod 9$, then there exists an overlattice of index 3 of $\langle 2d\rangle\oplus(E_6(-1)^{\oplus 3})'$ obtained by adding to the generators of $\langle 2d\rangle\oplus(E_6(-1)^{\oplus 3})'$ the class $L/3+\varepsilon$, where $\varepsilon$ is the generator of the discriminant form of $(E_6(-1)^{\oplus 3})'$. If $d\equiv 6\mod 9$ then the overlattice of index 3 of $\langle 2d\rangle\oplus(E_6(-1)^{\oplus 3})'$ is obtained by adding to the generators of $\langle 2d\rangle\oplus(E_6(-1)^{\oplus 3})'$ the class $L/3+2\varepsilon$. If $d\equiv 0 \mod 9$, there exists no proper overlattice of $\langle 2d\rangle\oplus (E_6(-1)^{\oplus 3})'$. The uniqueness of the overlattices constructed follows by observing that $\varepsilon$ is uniquely determined. 

The lattices $\langle 2d\rangle\oplus(E_6(-1)^{\oplus 3})'$ and $(\langle 2d\rangle\oplus (E_6(-1)^{\oplus 3})')'$ admit a primitive embedding in $\Lambda_{K3}$ and hence each of them determines a family of polarized K3 surfaces. The dimension of each of these families is 1, since the rank of the N\'eron--Severi group is 19, and each of these families corresponds to projective K3 surfaces which appear in a generalized Shioda--Inose structure of order 3. 
\endproof

\begin{corollary}\label{cor: SI order 3 implies Si order 2 }Let $X$ be a projective K3 surface admitting a generalized Shioda--Inose structure of order 3 $(A,\sigma_A,X,\sigma_X)$, then $X$ admits also a Shioda--Inose structure of order 2 $(B,X,\iota_B)$ and $T_B\simeq T_A$.
\end{corollary}
\proof

Let $X$ be a K3 surface which admits a generalized Shioda Inose structure of order 3. Then it admits a symplectic automorphism $\sigma_X$ and there exists an Abelian surface $A$ which admits an automorphism $\sigma_A$ such that $A/\sigma_A$ is birational to $X/\sigma_X$. Moreover $T_A\simeq T_X$.  Since $T_X\simeq T_A$, $T_X$ is primitively embedded in $U^{\oplus 3}\simeq H^2(A,\Z)$. So by Theorem \ref{theorem: recap 2 SI order 2}, $X$ admits a Shioda--Inose structure of order 2. I.e. $X$ admits a symplectic involution $\iota_X$ such that the minimal resolution of $X/\iota_X$ is a Kummer surface $Km_2(B)$ for a certain Abelian surface $B$ and $T_X\simeq T_B$. It follows that $T_A\simeq T_X\simeq T_B$.

\endproof

The Abelian surfaces $A$ and $B$  of the previous corollary are in the same family of polarized Abelian surfaces, since they have the same transcendental lattice. We are not using the same letter, because chosen a particular K3 surface $X$ nothing assures that the Abelian surfaces $A$ and $B$ are indeed isomorphic surfaces.

In view of Corollary \ref{cor: SI order 3 implies Si order 2 } one could expect that the lattice which characterizes the generalized Shioda--Inose structure of order 3 (i.e. the lattice $(E_6(-1)^{\oplus 3})'$) contains lattice which characterizes the Shioda--Inose structure of order 2 (i.e. the lattice $(E_8(-1)^{\oplus 2})'$), however the following result holds:
\begin{proposition}\label{prop: no abstract embedding E8E8 in E6E6E6}
	There exist no primitive embeddings of $E_8(-1)\oplus E_8(-1)$ in $(E_6(-1)\oplus E_6(-1)\oplus E_6(-1))'$.\end{proposition}
\proof If $E_8(-1)\oplus E_8(-1)$ were primitively embedded in $(E_6(-1)\oplus E_6(-1)\oplus E_6(-1))'$ we would denote by $R$ its orthogonal complement. Since $E_8(-1)$ is a unimodular lattice there are no overlattices of  
$E_8(-1) \oplus E_8(-1)\oplus R$ in which $R$ is primitively embedded. In particular, $E_8(-1)\oplus E_8(-1)\oplus R$ would be isometric to $(E_6(-1)\oplus E_6(-1)\oplus E_6(-1))'$. In particular, the number of roots of these two lattices have to be the same. As in proof of Proposition \ref{prop: negative definite lattices and genus}, a standard computation (using e.g. Magma, \cite{Ma}) allows one to compute the number of roots of $(E_6(-1)\oplus E_6(-1)\oplus E_6(-1))'$, which turns out to be the same of $E_6(-1)\oplus E_6(-1)\oplus E_6(-1)$, i.e. $3\cdot 72=216$.
On the other hand the roots of $E_8(-1)\oplus E_8(-1)\oplus R$ can not be less than the ones of $E_8(-1)\oplus E_8(-1)$ and $E_8(-1)\oplus E_8(-1)$ has $2\cdot 240=480$ roots.  Hence there exists no embedding of $E_8(-1)\oplus E_8(-1)$ in $(E_6(-1)\oplus E_6(-1) \oplus E_6(-1))'$. Again the unimodularity  of $E_8(-1)\oplus E_8(-1)$ implies that each embedding of   $E_8(-1)\oplus E_8(-1)$ in $(E_6(-1)\oplus E_6(-1) \oplus E_6(-1))'$ is necessarily primitive.\endproof

The previous result essentially says that there is no a canonical way to define a symplectic involution related to a Shioda--Inose structure (of order 2) on K3 surfaces which admits a symplectic automorphism of order 3 related with generalized Shioda Inose structure of order 3. Indeed the involution can not be defined on the lattice which characterizes the presence of the automorphism of order 3.

\subsubsection{Remarks on the case of order 6}\label{subsubsection: remarks n=6} 
We discussed in Remark \ref{rem: generalization definition of SI order n} a generalization of the definition of the Shioda--Inose structure to the order $n$ case, where $n$ is not necessarily 2 or 3.

The following result suggests that the generalization proposed is not the ``right" one, as discussed in Remark \ref{rem: no generalization to order 6 case}. To state the next result we use as definition of generalized Shioda--Inose structure of order 6 the one given in Remark \ref{rem: generalization definition of SI order n}.
\begin{prop}\label{prop: example order 6 fails}
Let $A$ be an Abelian surface with transcendental lattice $T_A\simeq U\oplus\langle 2\rangle$. Then there exists $Km_6(A)$, but $A$ is not a term of a generalized Shioda--Inose structure of order 6.
\end{prop}
\proof If $T_A\simeq A$, then $T_A$ can be primitively embedded in $U\oplus A_2$ and so $A$ admits an order 3 automorphism $\sigma_A$ such that $A/\sigma_A$ is birational to $Km_3(A)$. By Remark \ref{rem: autom order 6 on A}, it follows that $A$ admits also an automorphism of order 6, $\tau_A$, such that $A/\tau_A$ is birational to the K3 surface $Km_6(A)$. If there were a generalized Shioda--Inose structure of order 6 containing $A$, there would exist $X$ and $\tau_X$ such that $X$ is a K3 surface with $T_X\simeq T_A$ and $\tau_X$ is an order 6 symplectic automorphism of $X$ such that $X/\tau_X$ is birational to $Km_6(A)$. We claim that if $X$ is a K3 surface with $T_X\simeq T_A$, then $X$ does not admit symplectic automorphisms of order 6 and so it is not possible to find $\tau_X$ with the required properties and hence to construct the generalized Shioda--Inose structure of order 6.
Indeed if a K3 surface admits a symplectic automorphism of order 6, then its transcendental lattice admits a primitive embedding in $U^{\oplus 2}\oplus U(6)^{\oplus 2}$ (see e.g. \cite{GSnonprime}). But $T_X\simeq T_A\simeq U\oplus \langle 2\rangle$ can not be primitively embedded in $U^{\oplus 2}\oplus U(6)^{\oplus 2}$, since its orthogonal complement $R$ should be a lattice of rank 3 such that $A_{T_X\oplus R}=A_{T_X}\oplus A_{R}$ contains $(\Z/6\Z)^4$, and $A_{T_X}\simeq \Z/2\Z$, see also \cite[Proposition 1.15.1]{Nforms}.\endproof

\begin{rem}\label{rem: no generalization to order 6 case}{\rm Let $A$ be any Abelian surface for which exists the (generalized) Kummer surface $Km_n(A)$ of order $n$. If $n=2,3$, then there exists a (generalized) Shioda--Inose structure of order $n$ containing $A$ (see Theorems \ref{theorem: recap 2 SI order 2} and \ref{theorem: SI 2 order 3}). If there exists a notion of generalized Shioda--Inose structure of order $6$ such that for each Abelian surface for which exists the generalized Kummer surface $Km_6(A)$, there exists a generalized Shioda--Inose structure of order $6$ containing $A$, then this notion is not the one given in Remark \ref{rem: generalization definition of SI order n}. Indeed the surface $A$ considered in Proposition \ref{prop: example order 6 fails} admits a generalized Kummer construction of order 6, but not a generalized Shioda--Inose structure of order 6 in the sense of Remark \ref{rem: generalization definition of SI order n}.}\end{rem}

The Proposition result is surprising also in view of Corollary \ref{cor: SI order 3 implies Si order 2 }. Indeed let $A$ be an Abelian surface which admits an order 3 automorphism $\sigma_A$ such that $A/\sigma_A$ is birational to $Km_3(A)$.  The surface $A$ admits also an order 6 automorphism, $\tau_A=\sigma_A\circ\iota_A$, such that $A/\tau_A$ is birational to the K3 surface $Km_6(A)$, see Remark \ref{rem: autom order 6 on A}. Let $(A,\sigma_A, X,\sigma_A)$ be a generalized Shioda--Inose structure of order 3, so that $X$ is a K3 surface admitting a symplectic order 3 automorphism $\sigma_X$ and $T_X\simeq T_A$. By Corollary \ref{cor: SI order 3 implies Si order 2 }, $X$ admits also a Shioda--Inose structure $(A,X,\iota_X)$ of order 2. Hence $X$ admits two symplectic automorphisms $\sigma_X$ and $\iota_X$ of order 3 and 2 respectively. One can suppose that their composition is an automorphism of order 6 which can be used to construct a generalized Shioda--Inose structure of order 6. By Proposition \ref{prop: example order 6 fails} this supposition is wrong. Indeed if $A$ is as in Proposition \ref{prop: example order 6 fails}, then we can construct $X$, $\sigma_X$, $\iota_X$ as above but we know that $\sigma_X$ and $\iota_X$ do not commute (otherwise $\sigma_X\circ\iota_X$ would be a symplectic automorphism of order 6 on $X$ and by the proof of the proposition $X$ does not admit symplectic automorphisms of order 6).

We observe that in all the examples we will present we obtain $\langle\iota_X,\sigma_X\rangle=\mathfrak{S}_3$, the symmetric group of order 6.

\section{An example with Picard number 20}\label{sec: example picard number 20}
In \cite{GP} the K3 surface $X_3$ with $T_{X_3}\simeq A_2$ was considered as example of K3 surface admitting an order 3 symplectic automorphism. It is in particular observed that such an automorphism acts permuting cyclically three copies of $E_6(-1)$ contained in its N\'eron--Severi group. By Theorem \ref{theorem: SI 2 order 3}
this implies that it admits a generalized Shioda--Inose structure of order 3. Moreover, by Corollary \ref{cor: SI order 3 implies Si order 2 } $X_3$ also admits an involution $\iota_X$ attached to a Shioda--Inose structure (of order 2). In this section we describe explicitly both a generalized Shioda--Inose structure of order 3 and a Shioda--Inose structure of order 2 and we will prove the following results.

\begin{theorem}\label{theo: 3IV*}
	Let $X_3$ be the K3 surface with $T_{X_3}\simeq A_2$. Then $X_3$ admits an elliptic fibration $\mathcal{E}$ whose Weierstrass equation is \begin{equation}\label{eq: Weierstrass X3}y^2=x^3+(t^2-1)^4.\end{equation}
	\begin{enumerate}
		\item 
		The elliptic fibration admits a section of order 3 
		and the translation by this 3-torsion section is an order 3 symplectic automorphism $\sigma_X$. The desingularization of the quotient $X_3/\sigma_3$ is isomorphic to $X_3$.
		\item	Moreover $X_3$ admits a symplectic involution $\iota_X$ which acts on the equation \eqref{eq: Weierstrass X3} as $(x,y,t)\mapsto (x,-y,-t)$.
		\item The datum $(E_{\zeta}\times E_{\zeta}, \sigma_{E_{\zeta}\times E_{\zeta}}, X_3,\sigma_X)$ is a generalized Shioda--Inose structure of order 3 (see Section \ref{subsec: abelian Ezeta Ezeta and autom} for the definitions of the elements of the data).
		\item The datum  $(E_{\zeta}\times E_{\zeta} , X_3,\iota_X)$ is a Shioda-Inose structure (of order 2).	
		\item The group generated by $\iota_X$ and $\sigma_X$ is $\mathfrak{S}_3$.
		\item The quotient of $X_3$ by the group $\langle\iota_X,\sigma_X\rangle$ is isomorphic to $Km_2(E_\zeta\times E_\zeta)$.
	\end{enumerate}
\end{theorem}
\begin{corollary}
	With the same notation of Theorem \ref{theo: 3IV*}, one has $$A_2\simeq T_{E_\zeta\times E_{\zeta}}\simeq T_{X_3}\simeq T_{Km_3(E_\zeta\times E_\zeta)}\mbox{ and }A_2(2)\simeq T_{Km_2(E_\zeta\times E_\zeta)}.$$ 
\end{corollary}

\subsection{The Abelian surface $E_{\zeta}\times E_{\zeta}$ and K3 surface $X_3$}\label{subsec: abelian Ezeta Ezeta and autom}
The surface $X_3$ is well known: it is considered in \cite{SI} where the first construction of Shioda--Inose structure are provided;
it is studied in \cite{V} as one of the ``two most algebraic" K3 surfaces; the elliptic fibrations on such a surface were classified in \cite{Nis}.

Here we briefly recall its construction and its main properties.

Let $\zeta$ be a primitive third root of unity and $E_{\zeta}$ be the elliptic curve admitting an endomorphism, $\alpha_E$, of order 3, i.e. the lattice which defines $E_{\zeta}$ is $\langle 1,\zeta\rangle.$ A Weierstrass equation of $E_\zeta$ is $v^2=u^3-1$ (which exhibits such a curve as double cover of $\mathbb{P}^1_{u}$ and as a triple cover of $\mathbb{P}^1_{v}$).

Let us consider the Abelian surface $E_{\zeta}\times E_{\zeta}$: its N\'eron Severi group is generated by the classes of the following four curves:
$$E_1:=0_E\times E_{\zeta}, \ E_2:=E_{\zeta}\times 0_E,\ \Delta=\{(p,p), p\in E_{\zeta}\}, \ \Gamma=\{(p,\alpha_E(p)),\ p\in E_{\zeta}\}.$$ 
The intersection properties are the following $E_1\Delta=E_2\Delta=1$, $E_1\Gamma=E_2\Gamma=1$, $\Delta\Gamma=3$.
So $NS(E_{\zeta}\times E_{\zeta})=\langle E_1,E_2,\Delta, \Gamma\rangle\simeq U\oplus A_2(-1)$, where the isometry between the lattices are given by the choice of the new basis $\{E_1,E_2, E_1+E_2-\Delta, E_1+E_2-\Gamma\}$.
It follows $T_{E_{\zeta}\times E_{\zeta}}\simeq A_2$.

Moreover $E_{\zeta}\times E_{\zeta}$ admits an automorphism $\sigma_{E_{\zeta\times\zeta}}=\alpha_E\times\alpha_E^2$.

Let $X_3$ be the minimal resolution of the quotient $(E_{\zeta}\times E_{\zeta})/\langle\alpha_E\times\alpha_E^2\rangle$. It is proved in \cite{SI} that $X_3$ is a K3 surface whose transcendental lattice is $T_{X_3}\simeq A_2$ and whose N\'eron--Severi group is isometric to $U\oplus E_8(-1)\oplus E_8(-1)\oplus A_2(-1)$.
\subsection{The elliptic fibration $\mathcal{E}_3$ and the order 3 automorphism}

The K3 surface $X_3$ admits an elliptic fibration induced by the projection of $E_{\zeta}\times E_{\zeta}$ on one factor $E_{\zeta}$.

We denote  $\mathcal{E}_3:X_3\ra \mathbb{P}^1$ this elliptic fibration. By \cite{SI} (see also \cite{GP}) it is known that:
\begin{enumerate}
\item the reducible fibers of $\mathcal{E}_3$ are 3 fibers of type $IV^*$ (whose dual diagram is $\widetilde{E_6}$) and there are no other singular fibers;
\item the Mordell--Weil group, $MW(\mathcal{E}_3)$, of $\mathcal{E}_3$ is $\Z/3\Z$ and in particular there is a 3-torsion section;
\item the translation by the 3-torsion section is an order 3 symplectic automorphism $\sigma_X$ on $X_3$, which permutes cyclically three copies of $E_6(-1)$.
\end{enumerate}
To be more precise about the third point, we need to identify a set of generators of $NS(X_3)$: we denote $C_i^{(j)}$, $i=0,\ldots, 6$, $j=1,2,3$ the irreducible components of the $j$-th fiber of type $IV^*$, by $\mathcal{O}$ the zero section of the fibration and by $T_h$, $h=1,2$ the order three sections.
The intersection properties are the ones shown in the figure (if two curves are connected by a line, they intersect with multiplicity 1, otherwise they do not intersect, the self intersection of each class represented is $-2$):

\begin{align}\label{eq: fibers IV*} \xymatrix{\mathcal{O}\ar@{-}[r]&C_0^{(j)}\ar@{-}[r]&C_1^{(j)}\ar@{-}[r]&C_2^{(j)}\ar@{-}[r]\ar@{-}[d]&C_3^{(j)}\ar@{-}[r]&C_4^{(j)}\ar@{-}[r]&T_1\\
&&&C_5^{(j)}\ar@{-}[d]\\
&&&C_6^{(j)}\ar@{-}[d]\\&&&T_2}\end{align}

The automorphism $\sigma_X$ preserves each fiber and cyclically permutes the sections, so its action on the generators of the N\'eron--Severi group is $$S\mapsto T_1\mapsto T_2,\ \ C_0^{(j)}\mapsto C_4^{(j)}\mapsto C_6^{(j)},\ \ C_1^{(j)}\mapsto C_3^{(j)}\mapsto C_5^{(j)},\ \ C_2^{(j)}\mapsto C_2^{(j)}\mapsto C_2^{(j)},\ j=1,2,3.$$

The three copies of $E_6(-1)$ cyclically permuted by $\sigma_X$ are the following \begin{align}\label{eq: 3E6}\xymatrix{C_1^{(1)}\ar@{-}[r]&C_0^{(1)}\ar@{-}[r]&\mathcal{O}\ar@{-}[r]\ar@{-}[d]&C_0^{(2)}\ar@{-}[r]&C_1^{(2)}\\&&C_0^{(3)}\\
	C_3^{(1)}\ar@{-}[r]&C_4^{(1)}\ar@{-}[r]&T_1\ar@{-}[r]\ar@{-}[d]&C_4^{(2)}\ar@{-}[r]&C_3^{(2)}\\&&C_4^{(3)}\\C_5^{(1)}\ar@{-}[r]&C_6^{(1)}\ar@{-}[r]&T_2\ar@{-}[r]\ar@{-}[d]&C_6^{(2)}\ar@{-}[r]&C_5^{(2)}\\&&C_6^{(3)}}\end{align}

The quotient $X_3/\sigma_X$ is considered in \cite{GP}: it is a singular surface with 6 singular points of type $A_2$ and its minimal resolution $Y_3$ is a K3 surface with the same N\'eron--Severi group and transcendental lattice as $X_3$. Indeed it is isomorphic to $X_3$. We constructed $X_3$, and hence $Y_3$, as resolution of the quotient of the Abelian surface $E_\zeta\times E_\zeta$ by the order 3 automorphism  $\sigma_{E_\zeta\times E_\zeta}=\alpha\times \alpha^2$. So $Y_3\simeq Km_3(E_\zeta\times E_\zeta)$. We already observed that $T_{Y_3}\simeq T_{X_3}\simeq A_2\simeq T_{E_{\zeta}\times E_{\zeta}}$.

We conclude that $$(E_{\zeta}\times E_{\zeta},\sigma_{E_{\zeta}\times E_{\zeta}}, X_3,\sigma_X)\mbox{ is a generalized Shioda Inose structure of order 3 }$$ 
and the associated generalized Kummer surface $Km_3(E_{\zeta}\times E_{\zeta})\simeq Y_3\simeq X_3$.

\subsection{The involution $\iota_X$}
By the Corollary \ref{cor: SI order 3 implies Si order 2 } $X_3$ admits a Shioda--Inose structure (of order 2). Hence by the Theorem \ref{theorem: recap 2 SI order 2} there exists an involution $\iota_{X_3}$ which switches two orthogonal copies of $E_8(-1)$ in $NS(X_3)$.

Here we identify two orthogonal copies of $E_8(-1)$:

$$\xymatrix{C_3^{(3)}\ar@{-}[r]&C_4^{(3)}\ar@{-}[r]&T_1\ar@{-}[r]\ar@{-}[d]&C_4^{(2)}\ar@{-}[r]&C_3^{(2)}\ar@{-}[r]&C_2^{(2)}\ar@{-}[r]&C_1^{(2)}\\&&C_4^{(1)}\\
C_5^{(3)}\ar@{-}[r]&C_6^{(3)}\ar@{-}[r]&T_2\ar@{-}[r]\ar@{-}[d]&C_6^{(1)}\ar@{-}[r]&C_5^{(1)}\ar@{-}[r]&C_2^{(1)}\ar@{-}[r]&C_1^{(1)}\\&&C_6^{(2)}\\
}$$
By Theorems \ref{theorem: SI 2 order 3} there exists a symplectic involution $\iota_X$ on $X_3$ such that $(E_{\zeta}\times E_{\zeta}, X_3, \iota_X)$ is a Shioda--Inose structure. By Theorem \ref{theorem: SI 1 order 3}, $\iota_X^*$ permutes the two copies of $E_8(-1)$ and preserves the classes of $NS(X_3)$ which are orthogonal to these copies of $E_8(-1)$; so one knows the action of $\iota_X^*$ on the generators of the N\'eron--Severi of $X_3$, and hence on $NS(X_3)$. Since the automorphisms $\iota_X^*$ and $\sigma_X^*$ are symplectic they act trivially on the transcendental lattice, so one can describe the action of $\iota_X^*$ and $\sigma_X^*$ on $H^2(X_3,\Z)$. One can directly check that $$\langle \iota_X^*,\sigma_X^*\rangle= \mathfrak{S}_3\subset O(H^2(X_3,\Z))\mbox{ which implies that }\langle \iota_X,\sigma_X\rangle= \mathfrak{S}_3\subset Aut(X_3).$$

To give a more geometric description of $\iota_X$, we recall that $\mathcal{E}_3:X_3\ra\mathbb{P}^1_t$ has three fibers of type $IV^*$. We can assume that two of these three fibers are over $t=\pm 1$ and the third is over $t=\infty$. As a consequence, there exists an involution of the basis ($\mathbb{P}^1_t$) which permutes two of the reducible fibers and preserves the third, in particular, one can assume that this involution acts on $t$ as $t\mapsto -t$. The automorphisms which act only on the basis of an elliptic fibration are non-symplectic. Vice versa, the composition of the involution $t\mapsto -t$ on the basis and the elliptic involution on each fiber is a symplectic involution and we will show that it coincides with the involution $\iota_X$. 
To be more explicit one can write the Weierstrass equation of an elliptic K3 surfaces with 3 fibers of type $IV^*$ over the points $\pm 1$, $\infty$ (by using the information of \cite[Table IV.3.1]{Mi}) and one obtains:
$$y^2=x^3+(t^2-1)^4.$$
The symplectic form is $(dx\wedge dt)/y$ and the involution $(x,y,t)\mapsto (x,-y,-t)$ is symplectic.

It preserves the zero section of the fibration (because the elliptic involution on each elliptic curve preserves the zero of the group law), and it maps the section $T_1$ to the section $T_2$ (because maps each point of the elliptic fibers to its opposite).

Moreover, it maps components of the third reducible fiber to components of the third reducible fiber, in particular it preserves $C_0^{(3)}$, $C_1^{(3)}$, $C_2^{(3)}$ and switches $C_4^{(3)}$ with $C_6^{(3)}$ and $C_3^{(3)}$ with $C_5^{(5)}$ (here one is using that the action of an automorphism preserves the intersection properties).

Similarly one obtains the action on the components of the first two fibers, which are permuted. Hence one obtains $$C_0^{(1)}\leftrightarrow C_0^{(2)},\ C_1^{(1)}\leftrightarrow C_1^{(2)},\ C_2^{(1)}\leftrightarrow C_2^{(2)},\ C_3^{(1)}\leftrightarrow C_5^{(2)},\ C_4^{(1)}\leftrightarrow C_6^{(2)},\ C_5^{(1)}\leftrightarrow C_3^{(2)},\ C_6^{(1)}\leftrightarrow C_4^{(2)}.$$ 
So the involution $\iota_X$ coincides with the involution $(x,y,t)\mapsto (x,-y,-t)$.

\subsection{The equations of the quotient elliptic fibrations}
The quotient of $X_3$ by $\iota_X$ can be computed by the explicit equation of the Weierstrass form of $\mathcal{E}_3$ and of $\iota_X$.
By generalizing \cite[Proposition 5.1]{G} we have
\begin{proposition}\label{prop: eqution quotient involution}
Let $\mathcal{E}:X\ra\mathbb{P}^1_t$ be an elliptic fibration whose Weierstrass equation is $y^2=x^3+C(t)x^2+A(t)x+B(t)$ which admits the automorphism $\iota:(x,y,t)\mapsto (x,-y,-t)$.
Then: $$C(t)=c_4t^4+c_2t^2+c_0,\ \ A(t)=a_8t^8+a_6t^6+a_4t^4+a_2t^2+a_0,$$ $$B(t)=b_{12}t^{12}+b_{10}t^{10}+b_8t^8+b_6t^6+b_4t^4+b_2t^2+b_0;$$
the surface $X/\iota$ admits an elliptic fibration and a Weierstrass equation of a birational model of $X/\iota$ is
$$Y^2=X^3+C'(\tau)X^2+A'(\tau)X+B'(\tau)\mbox{ where }
C'(\tau)=c_4\tau^3+c_2\tau^2+c_0\tau,\ $$
 $$A'(\tau)=a_8\tau^6+a_6\tau^5+a_4\tau^4+a_2\tau^3+a_0\tau^2,$$ $$ B'(\tau)=b_{12}\tau^9+b_{10}\tau^8+b_8\tau^7+b_6\tau^6+b_4\tau^5+b_2\tau^4+b_0\tau^3.$$
\end{proposition}
\proof The proof is analogue to the one of \cite[Proposition 5.1]{G}: one first multiplies the equation of the elliptic fibration on $X$  by $t^6$ and then one uses the following invariant coordinates: $Y:=yt^3$, $X:=xt^2$, $\tau:=t^2$.\endproof 
Applying this proposition to our situation, i.e. to $A(t)=0$, $B(t)=(t^2-1)^4$ one obtains that a Weierstrass equation of a birational model of  $X_3/\iota_X$ is 
$$y^2=x^3+\tau^3(\tau-1)^4.$$
Let us denote by $Z$ the minimal resolution of $X_3/\iota_X$ and by $\mathcal{E}_Z$ the elliptic fibration on $Z$ whose equation is $y^2=x^3+\tau^3(\tau-1)^4$. The singular fibers of $\mathcal{E}_Z$ are $II^*+I_0^*+IV^*$. There are no sections of infinite order, since the rank of the trivial lattice is 20. There is a fiber of type $II^*$, so there can not be torsion sections. We conclude that the Mordell--Weil group is trivial and $NS(Z)\simeq U\oplus E_8(-1)\oplus D_4(-1)\oplus E_6(-1)$. The discriminant group is $\Z/2\times \Z/6\Z$ and the transcendental lattice is isometric to $A_2(2)$, see \cite[Table 2 line 323]{SZ}.
 
Since $T_Z\simeq A_2(2)\simeq T_{E_\zeta\times E_\zeta}(2)$ we have another proof of the fact that $Z\simeq Km(E_\zeta\times E_\zeta)$ and the involution $\iota_X(x,y,t)=(x,-y,-t)$ is the involution of a Shioda--Inose structure.

By the equation of the elliptic fibration $\mathcal{E}_X$, we compute an equation for $X_3/\sigma_X$ by applying the following well known result, see e.g. \cite[Section 3]{To}.
\begin{proposition}\label{prop: equation quotient trnaslation order 3}
	Let $E$ be an elliptic curve admitting a 3-torsion rational point. Then its equation can be written as $y^2=x^3+d^2$ or as $y^2=x^3+a^2(x-b)^2$ according to the fact that it admits or not an endomorphism of order 3. The 3-torsion point is $(x,y)=(0,d)$, $(x,y)=(0,ab)$ respectively. The quotient of $E$ by the translation by the 3-torsion point is $$y^2=x^3-27d^2, \ \ \ y^2=x^3-27a^2(x-4a^2-27b)^2$$
 respectively.
 \end{proposition} 
 
 In particular, one obtains that the equation of a birational model of $X_3/\sigma_X$ is \begin{equation}\label{eq: quotient}y^2=x^3-27(\tau^2-1)^4.
\end{equation}
One observes that the fibration described by the previous equation has three fibers of type $IV^*$, as expected. The involution $\iota_X$ induces the involution $(x,y,\tau)\mapsto(-x,y,-\tau)$ on the quotient  surface $X_3/\sigma_X$. This agrees with the fact that $\langle \sigma_X,\iota_X\rangle\simeq \mathfrak{S}_3$. Indeed $\langle \sigma_X\rangle$ is a normal subgroup of $\langle \sigma_X,\iota_X\rangle\simeq \mathfrak{S}_3$ and the group $\langle \sigma_X,\iota_X\rangle/\langle \sigma_X\rangle\simeq \Z/2\Z$ acts on $X_3/\sigma_X$. We denote $\overline{\iota_X}$ its generator, i.e. $\overline{\iota_X}$ is the involution induced by $\iota_X$ on the quotient $X_3/\sigma_X$.

We observe that all these automorphisms are symplectic and extend to automorphisms of the minimal resolution of the singular surface obtained as quotient.

To compute the equation of the quotient of the elliptic fibration \eqref{eq: quotient} by $\overline{\iota_X}$ one applies Proposition \ref{prop: eqution quotient involution} and one obtains $y^2=x^3-27T^3(T-1)^4$, which is an elliptic fibration over $\mathbb{P}^1_T$ with reducible fibers $II^*+IV^*+I_0^*$ and trivial Mordell--Weil lattice. Hence it is isomorphic to $Km_2(E_\zeta\times E_\zeta)$ which is the unique K3 surface admitting an elliptic fibration with these properties.

\subsection{The diagrams of the $\mathfrak{S}_3$-quotient and a final remark}
So one finds the following diagrams: on the left side the covers between the quotient surfaces are shown; the corresponding rational covers between the smooth models of the surfaces are on the right. In the last diagram the equations of a Weierstrass model of each of these surfaces are provided
$$
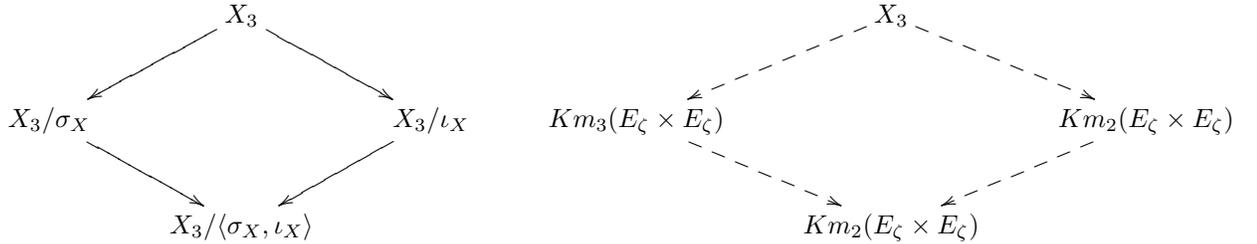
\xymatrix{&X_3\ar[dl]\ar[dr]&&&X_3\ar@{-->}[dl]\ar@{-->}[dr]\\
X_3/\sigma_X\ar[dr]&&X_3/\iota_X\ar[dl]&Km_3(E_\zeta\times E_\zeta)\ar@{-->}[dr]&&Km_2(E_\zeta\times E_\zeta)\ar@{-->}[dl]	\\
	&X_3/\langle\sigma_X,\iota_X\rangle&&&Km_2(E_\zeta\times E_\zeta)}$$
$$\xymatrix{
	&y^2=x^3+(t^2-1)^4\ar@{-->}[dl]\ar@{-->}[dr]\\
y^2=x^3-27(t^2-1)^4\ar@{-->}[dr]&&y^2=x^3+\tau^3(\tau-1)^4\ar@{-->}[dl]\\
&y^2=x^3-27T^3(T-1)^4}$$

We observe that the $3:1$ map $X_3/\iota\ra X_3/\langle \sigma_X,\iota_X\rangle$ is not Galois, all the other covers appearing in the diagram are cyclic (and in particular Galois).

\begin{rem}\label{rem: quotients different groups isomorphic}{\rm As observed in Remark \ref{rem: autom order 6 on A} $A:=E_{\zeta}\times E_{\zeta}$  admits an order 6 automorphism $\tau_A$ such that $A/\tau_A$ is birational to the K3 surface $Km_6(A)$. One has $\tau_A=\iota_A\circ\sigma_A$. We identified on the K3 surface $X_3$ two automorphisms $\sigma_X$ and $\iota_X$ which ``correspond" to $\sigma_A$ and $\iota_A$ by the Shioda--Inose construction. The composition $\sigma_X\circ\iota_X$ is not an automorphism of order 6 on $X_3$, in particular $(E_{\zeta}\times E_{\zeta},\iota_A\circ\sigma_A, X_3,\iota_X\circ \sigma_X)$ is not a generalized Shioda Inose of order 6, as defined in Remark \ref{rem: generalization definition of SI order n}, see also Section \ref{subsubsection: remarks n=6}. Nevertheless we claim that  $A/\langle \sigma_A,\iota_A\rangle$ is isomorphic to  $X_3/\langle\sigma_X,\iota_X\rangle$ (even if the groups $\Z/6\Z=\langle \sigma_A,\iota_A\rangle$ and $\mathfrak{S}_3=\langle \sigma_X,\iota_X\rangle$ are not isomorphic).
To prove this result, we denote $\iota_E$ the elliptic involution on $E_{\zeta}$ and, as above, $\alpha_E$ the order 3 endomorphism on $E_{\zeta}$, then we define an order 6 automorphism $\tau_E$ on $E_{\zeta}$ as $\tau_E=\sigma_E\circ\alpha_E$. The Abelian surface $E_{\zeta}\times E_{\zeta}$ admits an order 6 automorphism $\tau_A=\tau_E\times \tau_E^5$ which is the automorphism considered in Remark \ref{rem: autom order 6 on A}.
The projection $E_{\zeta}\times E_{\zeta}\ra E_{\zeta}$ induces a fibration $(E_{\zeta}\times E_{\zeta})/\tau_A\ra \mathbb{P}^1\simeq E_{\zeta}/\tau_E$. The surface $(E_{\zeta}\times E_{\zeta})/\tau_A=(E_{\zeta}\times E_{\zeta})/\langle \sigma_A,\iota_A\rangle$ has the following singularities: $A_5+4A_2+5A_1$ (which are the image of the points of $E_{\zeta}\times E_{\zeta}$ with non trivial stabilizer for $\tau_A$), see also \cite{Be}. With respect to the fibration  $(E_{\zeta}\times E_{\zeta})/\tau_A\ra \mathbb{P}^1$ one observes that there are three fibers which contain singular points: the first one contains a singularity of type $A_5$, one of type $A_2$, one of type $A_1$; the second one  contains 3 singularities of type $A_2$ and the last one contains 4 singularities of type $A_1$.  Hence the minimal resolution $Km_6(E_{\zeta}\times E_{\zeta})$ of $(E_{\zeta}\times E_{\zeta})/\tau_A$ has an elliptic fibration with 3 reducible fibers, one of type $II^*$, one of type $IV^*$ and one of type $I_0^*$. It follows that it is a K3 surface isomorphic to $Km_2(E_{\zeta}\times E_{\zeta})$ and in particular to $X_3/\langle \sigma_X,\iota_X\rangle$.}
\end{rem}
	
\section{Two examples with Picard number 19}\label{sec: examples with Picard number 19}
In the previous section we discussed in detail an example where the automorphism of order 3 involved in a generalized Shioda--Inose structure of order 3 is a translation by a 3-torsion section on an elliptic fibration. On the same elliptic fibration there is another automorphism, of order 2, which is the composition of an involution on the basis of the fibration with the elliptic involution on each fiber. We proved that the involution is a part of a Shioda--Inose structure (of order 2).

Here we consider several other examples of generalized Shioda--Inose structure of order 3 and Shioda-Inose structure of order 2, obtained in a similar way.
More precisely we consider K3 surfaces $X$ and an elliptic fibration $\mathcal{E}_X$ on each of these K3 surfaces such that $\MW(\mathcal{E})=\Z/3\Z$. So there exists an order 3 symplectic automorphism $\sigma_X$ which is the translation by a 3-torsion section. Moreover the K3 surface $X$ admits an order 2 automorphism, $\iota_X$, which preserves the fibration $\mathcal{E}$ and acts on the basis of the fibration as an involution and on the fiber as the elliptic involution. In this situation $\langle \sigma_X,\iota_X\rangle\simeq \mathfrak{S}_3$ and hence we have the following diagram:

$$
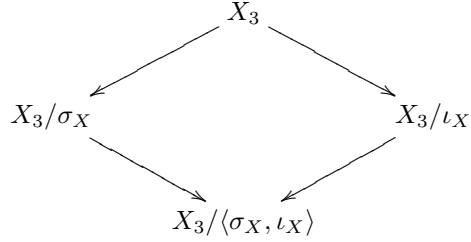
\xymatrix{&X_3\ar[dl]\ar[dr]&\\
	X_3/\sigma_X\ar[dr]&&X_3/\iota_X\ar[dl]\\
	&X_3/\langle\sigma_X,\iota_X\rangle}$$

In particular, we will show that in our examples the automorphisms $\sigma_X$ and $\iota_X$ are part of a generalized Shioda--Inose structure of order 3 and of a Shioda--Inose structure (of order 2) respectively. Hence there exists an Abelian surface $A$ such that the minimal resolution of the surfaces appearing in the previous diagram are as follows:
$$
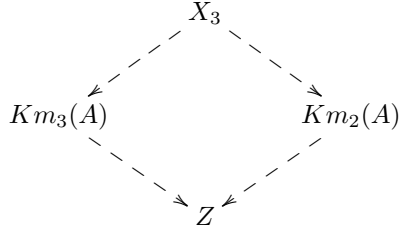
\xymatrix{&X_3\ar@{-->}[dl]\ar@{-->}[dr]\\
	Km_3(A)\ar@{-->}[dr]&&Km_2(A)\ar@{-->}[dl]	\\
	&Z}$$
 where $Z$ is a K3 surface birational to $X_3/\langle\sigma_X,\iota_X\rangle$.

\subsection{An elliptic fibration with $2IV^*+I_6$ as reducible fibers}\label{subsec: 2IV*I6}
In this section we consider a K3 surface $X$ whose transcendental lattice is $U\oplus \langle6\rangle$ and we prove the following theorem.
\begin{theorem}\label{theorem: 2IV*I6}
There exists a 1-dimensional family of K3 surfaces $X$ such that $NS(X)\simeq U\oplus E_8(-1)\oplus E_7(-1)\oplus A_2(-1)$ and $T_X\simeq U\oplus \langle 6\rangle\simeq A_2\oplus \langle -2\rangle$.
\begin{enumerate}
	\item Each K3 surface $X$ admits an elliptic fibration $\mathcal{E}$ with $2IV^*+I_6+2I_1$ as  singular fibers and $MW(\mathcal{E})=\Z/3\Z$. Its Weierstrass equation is 
	\begin{eqnarray}\label{eq: elliptic equation U+6}y^2=x^3+(t^2-1)^2(x-kt^4+kt^2)^2.\end{eqnarray}
	\item The translation by the 3-torsion section of $\mathcal{E}$ is an order 3 symplectic automorphism $\sigma_X$. The minimal resolution $Y$ of $X/\sigma_X$ admits an elliptic fibration whose Weierstrass equation is 	\begin{equation}\label{eq: 1 quotient U+6}y^2=x^3-27(t^2-1)^2(x-4(t^2-1)^2+27kt^4-27kt^2)^2.\end{equation}
	The reducible fibers are $2IV^*+2I_3+I_2$ and the Mordell--Weil group is $MW=\Z/3\Z$. Moreover, $NS(Y)\simeq U\oplus E_8(-1)\oplus E_6(-1)\oplus A_2(-1)\oplus A_1(-1)$ and $T_Y\simeq A_2\oplus\langle -6\rangle$. 
	\item The elliptic fibration $\mathcal{E}$ admits a symplectic involution $\iota_X:(x,y,t)\ra(x,-y,-t)$ and the minimal resolution of $X/\iota_X$ is a K3 surface admitting an elliptic fibration with $IV^*+I_3^*+I_0^*+I_1$ as singular fibers, trivial Mordell Weil lattice and whose equation is \begin{equation} y^2=x^3+t(t-1)^2x^2-2kt^3(t-1)^3x+k^2t^5(t-1)^4.\end{equation} In particular, the N\'eron--Severi group of this K3 surface is $U\oplus E_6(-1)\oplus D_7(-1)\oplus D_4(-1)$ and its transcendental lattice is $U(2)\oplus \langle 12\rangle\simeq (U\oplus\langle 6\rangle)(2)$. 

	\item There exists a 1-dimensional family of Abelian surfaces such that $T_A\simeq T_X$. Each Abelian surface $A$ in this family admits an order 3 automorphism $\sigma_A$ such that the resolution of $A/\sigma_A$ is a K3 surface $Km_3(A)$.
	\item $(A,\sigma_A,X,\sigma_X)$ is a generalized Shioda--Inose structure of order 3 and $(A, X, \iota)$ is a Shioda--Inose structure (of order 2). 
	\item The group $\langle \sigma_X,\iota_X\rangle$ is $\mathfrak{S}_3$ and the minimal model of the quotient surface $X/\langle \iota_X,\sigma_X\rangle$ is a K3 surface with an elliptic fibration with $IV^*+I_1^*+I_0^*+I_3$ as singular fiber, trivial Mordell Weil lattice and whose equation is \begin{equation}\label{eq: 2 quotient U+6} y^2=x^3-27t(t-1)^2x^2-54t^2(t-1)^3(27tk-4t+4)x-27t^3(t-1)^4(27tk-4t+4)^2.\end{equation}
    In particular, the N\'eron--Severi group of this K3 surface is $U\oplus E_6(-1)\oplus D_5(-1)\oplus D_4(-1)\oplus A_2(-1)$ and the transcendental lattice is  $A_2(2)\oplus \langle -12\rangle$.
\end{enumerate}
\end{theorem}
\proof $(1)$ The equation \eqref{eq: elliptic equation U+6} describes a K3 surface with an elliptic fibration whose singular fibers are $2IV^*+I_6+2I_1$ by \cite[Table IV.3.1]{Mi}. The fibration admits a 3-torsion section, whose equation is $t\mapsto(x(t),y(t))=(0,\pm(t^2-1)(x-kt^4+kt^2))$, so $\Z/3\Z\subset MW(\mathcal{E})$. Since the classes of the fiber, of the zero section and of the irreducible components of the reducible fibers span a lattice of rank 19, $\rho(X)\geq 19$ and $\rho(X)=19$ if and only if there are no sections of infinite order, by the Shioda--Tate formula, see e.g. \cite[Corollary VII.2.4]{Mi}. The equation \eqref{eq: elliptic equation U+6} depends on 1 parameter and the elliptic fibration is not isotrivial. Indeed for some specific values of $k$ (in particular $k=-\frac{4}{27}$), the fibration aquires a new singular fiber. This implies that the moduli space of the K3 surfaces admitting the elliptic fibration is 1-dimensional (if it were 0-dimensional, it is not possible that for certain values of the parameter the configuration of the reducible fibers changes).  Since the moduli space of the K3 surfaces $X$ has dimension 1, $\rho(X)= 19$ and there are no sections of infinite order. In particular, $MW(\mathcal{E})=\Z/3\Z$. The N\'eron--Severi of $X$ is generated by $F$, $\mathcal{O}$, $T_1$, $T_2$ and $C_i^{(j)}$ were $F$ is the class of a fiber, $\mathcal{O}$ is the class of the zero section, $T_i$ are the classes of the order 3 sections and $C_i^{(j)}$ are the irreducible components of the reducible fibers which do not intersect the zero section. One can now explicitly compute a basis of $NS(X)$ and its discriminant form, which suffices to determine both the N\'eron--Severi group and the transcendental lattice. The isometry between the two different lattices $U\oplus\langle 6\rangle$ and $A_2\oplus \langle -2\rangle$ follows by  \cite[Theorem 1.14.2]{Nforms}.

$(2)$ To compute the equation of the elliptic fibration on the minimal model of $X$, it suffices to apply Proposition \ref{prop: equation quotient trnaslation order 3}. Since there is a finite birational map between $X$ and the minimal resolution $Y$ of $X/\sigma_X$, it follows that $\rho(X)=\rho(Y)$ (see \cite{I}). By the equation of the elliptic fibration on $Y$, one determines the reducible fibers, by \cite[Table IV.3.1]{Mi}, and one deduces the presence of a 3-torsion section. Moreover, by the knowledge of the Picard number one deduces that there are no sections of infinite order and one deduces that the Mordell--Weil group of the fibration is $\Z/3\Z$. Again one has an explicit set of generators for $NS(Y)$, one extracts a basis, one finds the discriminant group and one deduces the description of the N\'eron--Severi group and of the transcendental lattice given in the statement.

$(3)$ The arguments are similar to the ones applied in point $(2)$, but in this case one has to use Proposition \ref{prop: eqution quotient involution} (instead of Proposition  \ref{prop: equation quotient trnaslation order 3}) in order to find the Weierstrass equation of the elliptic fibration.

$(4)$ The transcendental lattice of $X$ is $\langle -2\rangle \oplus A_2$, which can be embedded in $U\oplus A_2$ considering the basis $\{u_1-u_2,a_1,a_2\}$ of $U\oplus A_2$. Since $\langle -2\rangle \oplus A_2$ is primitively embedded in $U\oplus A_2\subset U^{\oplus 3}$, there exists an Abelian surface such that $T_A\simeq \langle -2\rangle \oplus A_2\simeq T_X$. Moreover, since $\langle -2\rangle \oplus A_2$ is primitively embedded in $U\oplus A_2$, by Theorem \ref{theorem: recap Kummer order 3}, $A$ admits an order 3 automorphism $\sigma_A$ such that $A/\sigma_A$ is birational to a K3 surface, denoted by  $Km_3(A)$. 
One can check that $T_{Km_3(A)}\simeq T_Y$, by applying the map $(\pi_A)_*$ to $T_A$. We recall that $(\pi_A)_*$ coincides with $\gamma$, see Definition \ref{def: gamma}. So $\gamma(\langle u_1-u_2,a_1,a_2\rangle)=\langle w_1-w_2,(b_1+2b_2)/3,b_2\rangle\simeq \langle -6\rangle \oplus A_2$.

$(5)$ We showed that $X$ admits an order 3 symplectic automorphism $\sigma_X$ such that the minimal resolution $Y$ of $X/\sigma_X$ has $T_Y\simeq A_2\oplus\langle -6\rangle$ and we also showed that there exists $A$ such that $T_A\simeq T_X$ and $T_{Km_3(A)}\simeq T_Y$. Hence $(A,\sigma_A,X,\sigma_X)$ is a generalized Shioda--Inose structure. Similarly, given $A$ such that $T_A\simeq T_X$, its Kummer surface $Km_2(A)$ has the same transcendental lattice as the one of the minimal resolution of $X/\iota_X$, i.e. $T_X(2)$ and then $(A,X,\iota_X)$ is a Shioda--Inose structure (of order 2). 

To make the construction of the (generalized) Shioda--Inose structures more explicit we show 3 copies of $E_6(-1)$ permuted by $\sigma_X^*$ and two copies of $E_8(-1)$ switched by $\iota_X^*$. The notation are as before, i.e. $C_i^{(j)}$, $i=0,\ldots,6$, $j=1,2$ are the components of the two reducible fibers of type $IV^*$ and their intersection properties are as in \eqref{eq:E6}. Moreover $C_i^{(3)}$, $i=0,\ldots, 5$ are the components of the fiber of type $I_6$, with the convention $C_i^{(3)}C_j^{(3)}=1$ if $|i-j|=1$, $i,j\in \Z/6\Z$. Hence $\sigma_X^*$ permutes cyclically the following copies of $E_6(-1)$:

\begin{align}\label{eq: 3E6 second case}\xymatrix{C_1^{(1)}\ar@{-}[r]&C_0^{(1)}\ar@{-}[r]&\mathcal{O}\ar@{-}[r]\ar@{-}[d]&C_0^{(2)}\ar@{-}[r]&C_1^{(2)}\\&&C_0^{(3)}\\
		C_3^{(1)}\ar@{-}[r]&C_4^{(1)}\ar@{-}[r]&T_1\ar@{-}[r]\ar@{-}[d]&C_4^{(2)}\ar@{-}[r]&C_3^{(2)}\\&&C_2^{(3)}\\C_5^{(1)}\ar@{-}[r]&C_6^{(1)}\ar@{-}[r]&T_2\ar@{-}[r]\ar@{-}[d]&C_6^{(2)}\ar@{-}[r]&C_5^{(2)}\\&&C_4^{(3)}}\end{align}
	
	The involution $\iota_X^*$ switches the following two  orthogonal pairs of $E_8(-1)$:
	
	$$
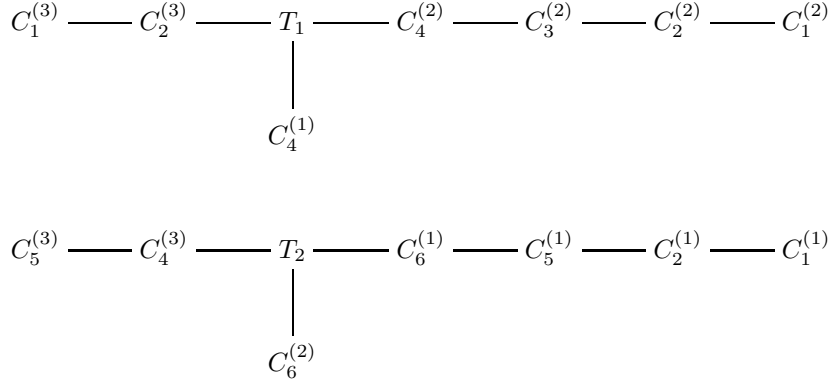
\xymatrix{C_1^{(3)}\ar@{-}[r]&C_2^{(3)}\ar@{-}[r]&T_1\ar@{-}[r]\ar@{-}[d]&C_4^{(2)}\ar@{-}[r]&C_3^{(2)}\ar@{-}[r]&C_2^{(2)}\ar@{-}[r]&C_1^{(2)}\\&&C_4^{(1)}\\
		C_5^{(3)}\ar@{-}[r]&C_4^{(3)}\ar@{-}[r]&T_2\ar@{-}[r]\ar@{-}[d]&C_6^{(1)}\ar@{-}[r]&C_5^{(1)}\ar@{-}[r]&C_2^{(1)}\ar@{-}[r]&C_1^{(1)}\\&&C_6^{(2)}\\
	}$$

$(6)$ Since there is an explicit equation of $X$ and of its automorphisms $\sigma_X$ and $\iota_X$, one observes that $\sigma_X^3=\iota_X^2=id$ and that $\sigma_X\iota_X=\iota_X\sigma_X^2$, so $\langle \sigma_X,\iota_X\rangle\simeq \mathfrak{S}_3$. To find $X/\langle \sigma_X,\iota_X\rangle$, one considers the iterated quotients $(X/\sigma_X)/\iota_X$. So it suffices to apply Proposition \ref{prop: eqution quotient involution} to the equation \eqref{eq: 1 quotient U+6} of the elliptic fibration on $Y$. One obtains the equation \eqref{eq: 2 quotient U+6} and one deduces the reducible fibers and as before one knows that the Picard number of this surface is 19 and so the Mordell--Weil group is of finite rank. Moreover there can not be torsion sections, because the fibers of type $IV^*$ allows only 3-torsion sections and the fibers of type $I_n^*$ do not admit 3-torsion section.
Hence one finds the N\'eron--Severi group and computes the  transcendental lattice. 
To describe the geometry of the involution induced by $\iota_X$ on $Y$ and to find the transcendental lattice of the quotient in a different way, we observe the following. The involution $\iota_X$ induces on $Y$ an involution, denoted by $\iota_Y$. This involution switches the two fibers of type $IV^*$ and the two fibers of type $I_3$. Denoting $T_1^Y$, $T_2^Y$ the torsion section of the elliptic fibration on $Y$ and by $D_i^{(j)}$ the irreducible components of $j$-th reducible fiber of $Y$, the non trivial intersections are the following:
$$\begin{array}{c}D_i^{(j)}D_{i+1}^{(j)}=1,\ j=1,2,\ i=0,1,2,3,5,\ \  D_2^{(j)}D_5^{(j)}=1,\ j=1,2\\
D_h^{(k)}D_{h+1}^{(k)}=1,\ h=0,1,2\ h\in\Z/3\Z,\ k=3,4,\ D_0^{(5)}D_1^{(5)}=2,\\ T_1^{Y}D_4^{(j)}=T_1^YD_1^{(k)}=T_1^YD_0^{(5)}=1,\ j=1,2,\ k=3,4\\ T_2^{Y}D_6^{(j)}=T_2^YD_2^{(k)}=T_2^YD_0^{(5)}=1,\ j=1,2,\ k=3,4.\end{array}$$
The involution $\iota_Y$ acts as follows:
$$\begin{array}{c}T_1^Y\leftrightarrow T_2^Y,\ D_0^{(1)}\leftrightarrow D_0^{(2)},\ D_1^{(1)}\leftrightarrow D_1^{(2)},\ D_2^{(1)}\leftrightarrow D_2^{(2)},\ D_3^{(1)}\leftrightarrow D_5^{(2)},\ D_4^{(1)}\leftrightarrow D_6^{(2)},\\ D_5^{(1)}\leftrightarrow D_3^{(2)},\ D_6^{(1)}\leftrightarrow D_4^{(2)},\ D_1^{(3)}\leftrightarrow D_2^{(4)},\ \ D_2^{(3)}\leftrightarrow D_1^{(4)},\ \ D_0^{(3)}\leftrightarrow D_0^{(4)}.\end{array}$$
Hence $\iota_Y$ switches the two following orthogonal copies of $E_8(-1)$:
	$$\xymatrix{D_0^{(3)}\ar@{-}[r]&D_1^{(3)}\ar@{-}[r]&T_1\ar@{-}[r]\ar@{-}[d]&D_4^{(2)}\ar@{-}[r]&D_3^{(2)}\ar@{-}[r]&D_2^{(2)}\ar@{-}[r]&D_1^{(2)}\\&&D_4^{(1)}\\
	D_0^{(4)}\ar@{-}[r]&D_2^{(4)}\ar@{-}[r]&T_2\ar@{-}[r]\ar@{-}[d]&D_6^{(1)}\ar@{-}[r]&D_5^{(1)}\ar@{-}[r]&D_2^{(1)}\ar@{-}[r]&D_1^{(1)}\\&&D_6^{(2)}\\
}$$

So the transcendental lattice of the minimal resolution $Z$ of $Y/\iota$ is such that $T_Z\simeq T_Y(2)\simeq A_2(2)\oplus \langle -12\rangle$, by Theorem \ref{theorem: recap 1 SI order 2}.

\subsection{An elliptic fibration with $I_{18}$ as reducible fiber}\label{subsec: I18}

\begin{theorem}
	There exists a 1-dimensional family of K3 surfaces $X$ such that $NS(X)\simeq U\oplus E_8(-1)\oplus E_8(-1)\oplus A_1(-1)$ and $T_X\simeq U\oplus \langle 2\rangle$.
	\begin{enumerate}
		\item Each K3 surface $X$ admits an elliptic fibration $\mathcal{E}$ with $I_{18}+6I_1$ as  singular fibers and $MW=\Z/3\Z$. 
		\item The translation by the 3-torsion section of $\mathcal{E}$ is an order 3 symplectic automorphism $\sigma_X$. The minimal resolution $Y$ of $X/\sigma_X$ admits an elliptic fibration whose singular fibers are $I_6+6I_3$ and the Mordell--Weil group is $MW=\Z/3\Z\times \Z/3\Z$. The transcendental lattice of this surface is $U(3)\oplus \langle 6\rangle\simeq (U\oplus \langle 2\rangle)(3)$;
		\item The elliptic fibration $\mathcal{E}$ admits a symplectic involution which is the composition of the elliptic involution on each fiber and an involution on the basis which preserves the point over which $\mathcal{E}$ has the reducible fiber of type $I_{18}$. The minimal resolution of $X/\iota_X$ is a K3 surface admitting an elliptic fibration with $I_9^*+I_0^*+3I_1$ as singular fibers and trivial Mordell--Weil lattice. The N\'eron--Severi group of this K3 surface is $U\oplus D_{13}(-1)\oplus D_4(-1)$ and its transcendental lattice is $U(2)\oplus \langle 4\rangle\simeq (U\oplus\langle 2\rangle)(2)$.
		\item There exists a 1-dimensional family of Abelian surfaces such that $T_A\simeq T_X$. Each Abelian surface $A$ in this family admits an order 3 automorphism such that the resolution of $A/\sigma_A$ is a K3 surface.
		\item $(A,\sigma_A,X,\sigma_X,)$ is a generalized Shioda--Inose structure of order 3 and $(A, X, \iota_X)$ is a Shioda--Inose structure (of order 2). 
		\item The group $\langle \sigma_X,\iota_X\rangle$ is $\mathfrak{S}_3$ and the minimal model of the quotient surface $X/\langle \sigma_X,\iota_X\rangle$ is a K3 surface with an elliptic fibration with $I_3^*+I_0^*+3I_3$ as singular fibers and trivial Mordell Weil lattice. The N\'eron--Severi group of this K3 surface is $U\oplus D_7(-1)\oplus D_4(-1)\oplus A_2(-1)^{\oplus 3}$ and the transcendental lattice is  $U(6)\oplus \langle 12\rangle=\left(U\oplus \langle 2\rangle\right)(6)$.
	\end{enumerate}
\end{theorem}
\proof The proof is similar to the ones of the Theorems \ref{theo: 3IV*} and \ref{theorem: 2IV*I6}.

Even without writing explicitly the Weierstrass equations of the elliptic fibration involved in the constructions, one is able to determine the type of singular fibers of the quotient elliptic fibration. Indeed if $\mathcal{E}$ is an elliptic fibration with a $n$ torsion section, the quotient of $\mathcal{E}$ by the translation by this torsion section is still an elliptic fibration and its reducible fibers can be deduced by the one of $\mathcal{E}$. We are interested in fibers of type $I_n$, so we can apply \cite[Lemma 2.2]{GSnonprime} to determine the fibers of $Y$. By \cite[Line 2942 Table 1]{Shim}, if the configuration of the singular fibers on an elliptic K3 surface is $I_6+6I_2$, then the torsion part of its Mordell--Weil group is $(\Z/3\Z)^2$. As above, we observe that $\rho(X)=\rho(Y)$, so there are no sections of infinite order. This allows to find a basis of the N\'eron--Severi group, to compute its discriminant and to deduce $T_X$ by this.

Similarly, $\iota_X$ identifies pairs of fibers and preserves two fibers, in particular the one of type $I_{18}$ and one smooth fiber $I_0$. It acts as the elliptic involution on each of them, hence the quotient fibers are $I_9^*$ and $I_0^*$.

Moreover, to compute the transcendental lattice, one checks from the discriminant form that the transcendental lattice is a multiple of an even lattice $L$, i.e. one obtains that the transcendental lattice $T$ is $T\simeq L(n)$. The length of $L$ is smaller than the one of $T$, and so one can apply the results in \cite{Nforms} to determine the unique lattice $L$ with the required discriminant group and form. We apply this argument to find the transcendental lattice of the minimal resolution of $X/\langle \sigma_X,\iota_X\rangle\simeq (X/\sigma_X)/\iota_X$.

The rest of the proof is analogous to the one of Theorem \ref{theorem: 2IV*I6}.

To be more explicit about the geometry of the (generalized) Shioda--Inose structures we show the three copies of $E_6(-1)$ cyclically permuted by $\sigma_X^*$ and the two copies of $E_8(-1)$ switched by $\iota_X^*$. Let us denote by $F$ the class of the fiber of the elliptic fibration $\mathcal{E}$, by $\mathcal{O}$ the zero section, by $T_1$, $T_2$ the two 3-torsion sections and by $C_i$, $i=0,\ldots 17$ the irreducible components of the reducible fiber.
The non trivial intersections between these curves are $\mathcal{O}^2=T_1^2=T_2^2=C_i^2=-2$, $C_iC_{i+1}=1$, $i=0,\ldots, 17$, $i\in\Z/18\Z$, $\mathcal{O}C_0=T_1C_6=T_2C_{12}=1$. 
The automorphism $\sigma_X^*$ acts as follows: 
$$\mathcal{O}\ra T_1\ra T_2,\ \ C_i\ra C_{i+6}\ra C_{i+12},\ \ i\in\Z/18\Z.$$
So the following three copies of $E_6(-1)$ are cyclically permuted by $\sigma_X^*$

\begin{align}\xymatrix{C_2\ar@{-}[r]&C_1\ar@{-}[r]&C_0\ar@{-}[r]\ar@{-}[d]&C_{17}\ar@{-}[r]&C_{16}\\&&\mathcal{O}\\
		C_8\ar@{-}[r]&C_7\ar@{-}[r]&C_6\ar@{-}[r]\ar@{-}[d]&C_5\ar@{-}[r]&C_4\\&&T_1\\C_{14}\ar@{-}[r]&C_{13}\ar@{-}[r]&C_{12}\ar@{-}[r]\ar@{-}[d]&C_{11}\ar@{-}[r]&C_{10}\\&&T_2}\end{align}

The involution $\iota_X^*$ acts as follows $$T_1\lra T_2,\ \ C_i\lra C_{18-i},\ \ i\in\Z/18\Z$$ and fixed $\mathcal{O}$, $C_0$ and $C_9$. Hence $\iota_X^*$
switches the following two orthogonal pairs of $E_8(-1)$:

$$\xymatrix{C_8\ar@{-}[r]&C_7\ar@{-}[r]&C_6\ar@{-}[r]\ar@{-}[d]&C_5\ar@{-}[r]&C_4\ar@{-}[r]&C_3\ar@{-}[r]&C_2\\&&T_1\\
	C_{10}\ar@{-}[r]&C_{11}\ar@{-}[r]&C_{12}\ar@{-}[r]\ar@{-}[d]&C_{13}\ar@{-}[r]&C_{14}\ar@{-}[r]&C_{15}\ar@{-}[r]&C_{16}\\&&T_2\\
}$$

\begin{remark}{\rm As already mentioned in the introduction, in \cite{Sc} the author considers infinite series of K3 surfaces $X$ which are part of Shioda--Inose structures (of order 2). In particular, $X$ admits an involution $\iota_X$ such that the minimal model of $X/\iota_X$ is isomorphic to $Km_2(A)$. The author also shows that  all the K3 surfaces $X$ constructed are in the middle of a sandwich made up of the Kummer surface $Km_2(A)$. This means that there is an involution $\iota_{Km}\in\Aut(Km_2(A))$ such that the minimal resolution of $Km_2(A)/\iota_{Km}$ is isomorphic to $X$, i.e. there exists the following diagram: $$Km_2(A)\dashrightarrow X\dashrightarrow Km_2(A)$$
		where the dash arrows are rational maps generically $2:1$.

All the examples considered in this paper have the same property, i.e. if $(A,\sigma_a,X,\sigma_X)$ is a generalized Shioda--Inose structure of order 3, then on $Km_3(A)$ there exists an order 3 automorphism $\sigma_{Km}$ such that $Km_3(A)/\sigma_{Km}$ is isomorphic to $X$, i.e. there exists the following diagram:
$$Km_3(A)\dashrightarrow X\dashrightarrow Km_3(A)$$
where the dash arrows are rational maps generically $3:1$.

This follows directly from the fact that the automorphism $\sigma_X$ is a translation by a 3-torsion section. Indeed in \cite{G} it is proved that if $S$ is a K3 surface with an elliptic fibration with an $n$ torsion section, denoted by $\sigma_S$ the translation by  this torsion section, the minimal resolution $W$ of $S/\sigma_n$ admits an elliptic fibration with an $n$-torsion section. Denoted by $\sigma_W$ the translation by this torsion section, the minimal resolution of $W/\sigma_W$ is isomorphic to $S$. This means that $S$ is in the middle of a sandwich made up of $W$.
} \end{remark}

\end{document}